\documentclass[11pt]{amsart}
\usepackage{amstext,amssymb,amsmath,amsbsy,color}

\usepackage{hyperref}
\usepackage{amscd}
\usepackage{amsfonts}
\usepackage{indentfirst}
\usepackage{verbatim}
\usepackage{amsmath}
\usepackage{amsthm}
\usepackage{enumerate}
\usepackage{graphicx}
\usepackage[OT1]{fontenc}
\usepackage[latin1]{inputenc}
\usepackage[english]{babel}
\usepackage{amssymb}

\newtheorem{theorem}{Theorem}
\newtheorem{lemma}{Lemma}

\newtheorem{proposition}{Proposition}
\newtheorem{corollary}{Corollary}
\newtheorem{remark}{Remark}

\numberwithin{equation}{section}

\newcommand{\proofend}{\hfill $\Box$ }
\newcommand{\mint}{\rule[1.1mm]{2.45mm}{.1mm}\hspace{-3.35mm}\int}

\newcommand{\supp}{\operatorname{supp}}
\newcommand{\dist}{\operatorname{dist}}

\newcommand{\dive}{\operatorname{div}}

\newcommand{\eps}{\varepsilon}

\newcommand{\loc}{_{loc}}

\newcommand{\mS}{\mathbb{S}}
\newcommand{\mR}{\mathbb{R}}

\newcommand{\mC}{\mathbb{C}}

\newcommand{\mc}{\mathrm{c}}

\title[Recovering potential from Cauchy data]{Recovering a potential from Cauchy data \\ via complex geometrical optics solutions}  

\author{Hoai-Minh Nguyen}
\thanks{Ecole Polytechnique F\'ed\'erale de Lausanne, EPFL SB MATHAA CAMA, Station 8,  CH-1015 Lausanne, Switzerland, {\tt{hoai-minh.nguyen@epfl.ch}}}
\author{Daniel Spirn}
\thanks{School of Mathematics,
University of Minnesota,
Minneapolis, MN 55455,
USA, \tt{spirn@math.umn.edu}}
\date{\today}

\begin{document}    

\begin{abstract}
This paper is devoted to  the problem of recovering a potential $q$ in a domain in $\mathbb{R}^d$ for $d \geq 3$ from the Dirichlet to Neumann map.  This problem is related to the inverse Calder\'on conductivity problem   via the Liouville transformation.  It is known from
the work of Haberman and Tataru  \cite{HabermanTataru11} and  Nachman and Lavine \cite{LavineNachman}  that  uniqueness holds for the class of conductivities of one derivative and  the class of $W^{2,d/2}$ conductivities respectively.  The proof of Haberman and Tataru is based on the construction of  complex geometrical optics (CGO) solutions initially suggested by  Sylvester and Uhlmann \cite{SylvesterUhlmann87}, in functional spaces introduced by Bourgain \cite{B}. The proof of the second result, in the work of Ferreira et al. \cite{FerreiraKenigSalo},  is based on the construction of  CGO solutions via Carleman estimates.  
 The main goal of the paper is to understand whether or not an approach which is based on the construction of CGO solutions in the spirit of Sylvester and Uhlmann and involves only standard Sobolev spaces can be used to obtain these results. In fact, 
 we are able to  obtain a new proof of uniqueness for the Calder\'on problem for 1) a slightly different class as the one in \cite{HabermanTataru11}, and for 2) the class of $W^{2,d/2}$ conductivities. The proof of statement 1) is based on a new estimate
for CGO solutions  and some averaging  estimates  in the same spirit as in \cite{HabermanTataru11}. The  proof of statement 2)  is on the one hand  based on a generalized Sobolev inequality due to Kenig et al. \cite{KRS} and on another hand, only involves standard estimates for CGO solutions   \cite{SylvesterUhlmann87}. We are also able to prove the uniqueness of a potential for  3) the class of $W^{s, 3/s}$ ($\supsetneqq W^{2, 3/2}$)  conductivities with $3/2 < s < 2$ in three dimensions. As far as we know, statement 3) is new. 

\end{abstract}
   
\maketitle 

\section{Introduction}

Let $\Omega$ be a bounded domain in $\mR^d$ $(d \ge 3)$ with $C^1$ boundary and let $q \in L^{d/2}(\Omega)$, an assumption that will be weaken later. We consider the Dirichlet to Neumann map  $\Lambda_{q}: H^{1/2}(\partial \Omega) \to H^{-1/2}(\partial \Omega)$ given by
\begin{equation*}\label{def-Lambda}
\Lambda_{q} (f) = g,
\end{equation*}
where
\begin{equation*}
g =  \frac{\partial v}{\partial \eta} \Big|_{\partial \Omega}, 
\end{equation*}
and  $v \in H^1(\Omega)$ is the unique solution to the system
\begin{equation*}\left\{
\begin{array}{cl}
\Delta v - q v  = 0 & \mbox{ in } \Omega, \\[6pt]
v = f & \mbox{ on } \partial \Omega.
\end{array}\right.
\end{equation*}
Here and in what follows $\eta$ denotes a unit normal vector directed into the exterior of $\Omega$. We assume here that 0 is not a Dirichlet eigenvalue for this problem; this implies $\Lambda_{q}$ is well-defined (this assumption is not essential and is discussed later in Remark~\ref{rem1}). In this paper, we are interested in the injectivity of $\Lambda_{q}$ for $d \ge 3$.  This problem has a connection to the inverse conductivity problem posed by Calder\'on in \cite{Calderon}. In \cite{Calderon} Calder\'on asked whether one can determine $\gamma \in L^{\infty}(\Omega)$ with $\mbox{essinf}_{\Omega} \gamma > 0$ from its Dirichlet to Neumann map $DtN_{\gamma}: H^{1/2}(\partial \Omega) \to H^{-1/2}(\partial \Omega)$ given by 
\begin{equation*}
DtN_{\gamma}(f) = \gamma \frac{\partial u}{\partial \eta}, 
\end{equation*}
where $u \in H^{1}(\Omega)$ is the unique solution to the equation 
\begin{equation*}
\dive(\gamma \nabla u) = 0 \mbox{ in } \Omega \mbox{ and } u = f \mbox{ on } \partial \Omega. 
\end{equation*}
In the same paper, Calder\'on proved the injectivity of the derivative  of the map $\gamma \to DtN_{\gamma}$ at $\gamma = $ constant. 
Kohn and Vogelius \cite{KohnVogelius84, KohnVogelius85}  showed that if $\partial \Omega $ is $C^\infty$ then $\Lambda_{q}$ determines $q$ and all its derivatives on $\partial \Omega$ and then used this to prove  uniqueness for the class of piecewise analytic coefficients.  Sylvester and Uhlmann \cite{SylvesterUhlmann87} proved that $\Lambda_{q}$ uniquely determines $q$ if $q \in C^\infty$; their method also gave the injectivity of $\Lambda_{q}$ for $q \in  L^{\infty}$ (see also \cite{Nachman88}). In \cite{SylvesterUhlmann87}, they introduced the concept of {\it complex geometrical optics} (CGO) solutions which plays an important role in establishing the uniqueness for inverse problems for $d \ge 3$. In one direction,  the $L^{\infty}$ uniqueness result was improved by  Chanillo, and Kenig and Jerrison in \cite{Chanillo}  and Lavine and Nachman in \cite{LavineNachman}. 
In \cite{Chanillo}, Chanillo established the injectivity of $\Lambda_{q}$ for $q \in L^{d/2}$ with small norm and (in the same paper) 
Kenig and Jerrison obtained the injectivity of $\Lambda_{q}$ for  $q \in L^{p}$ for any $p>d/2$.  In \cite{LavineNachman}, the authors announced the injectivity of  $\Lambda_{q}$ holds for $q \in L^{d/2}$. Recently, this result has been extend by  Ferreira et al. in \cite{FerreiraKenigSalo} for compact Riemannian manifolds with boundary which are conformally embedded in a product of the Euclidean line and a simple manifold. Their technique is  based on Carleman estimates. In another direction,  the injectivity of $\Lambda_{q}$ was  established for $q \in B^{-s}_{\infty, 2}$ \footnote{$B^s_{p, q}$ denotes the Besov spaces.} ($0 < s < 1/2$), $q \in B^{-1/2}_{\infty, 2}$, and for $q \in W^{-1/2, s}$ ($s > 2d$) by Brown in \cite{Brown96}, P\"aiv\"arinta et al. in \cite{PaivarintaPanchenkoUhlmann03}, and Brown and Torres in \cite{BT}, respectively. Recently, Haberman and Tataru in \cite{HabermanTataru11} established  the injectivity  of $DtN_{\gamma}$ (Calder\'on's problem) for $\gamma \in C^1(\Omega)$ or $\gamma \in W^{1, \infty}({\overline{\Omega}})$ with a smallness assumption on the derivative. The corresponding uniqueness result for $\Lambda_{q}$ would hold for $q \in W^{-1, \infty}$ 
with some kind of smallness assumption; however, obtaining this conclusion from their approach is not clear to us.  
The approach in \cite{Brown96,PaivarintaPanchenkoUhlmann03,BT} is via CGO solutions. The approach due to Haberman and Tataru is also via CGO solutions;  the novelty in  their approach stems from their  use of  weighted spaces and  averaging arguments. Some refinements for piecewise smooth potentials $q$ can be found in references therein (see also \cite{Isakov}). We note that the result of Lavine and Nachman is not a consequence of the one of Haberman and Tataru and vice versa since $L^{d/2}(\Omega) \not \in W^{-1, \infty}(\Omega)$ and $W^{-1, \infty}(\Omega) \not \in L^{d/2}(\Omega)$.  In  dimension $2$, the injectivity of $\Lambda_q$ was established by Astala and P\"aiv\"arinta in \cite{AstalaPaivarinta06}. Previous contributions in the $2d$ case can be found in \cite{Nachman96, BrownUhlmann97} and references therein.

\medskip
The standard method to establish  uniqueness for the Calder\'on problem is to prove the injectivity of  $\Lambda_{q}$. This can be done by the Liouville transform and using the fact that one can recover the boundary data from the Dirichlet to Neumann map since
\begin{equation}\label{eq1}
\Delta v - q v = 0 \mbox{ in } \Omega
\end{equation}
if and only if 
\begin{equation*}
\dive(\gamma \nabla u ) = 0 \mbox{ in } \Omega,
\end{equation*}
where $u = \gamma^{1/2} v$ and $q = \frac{\Delta \gamma^{1/2}}{\gamma^{1/2}}$. 
It is known that  (see e.g. \cite[(5.0.4)]{Isakov-book}) if 
\begin{equation*}
\Lambda_{q_{1}} = \Lambda_{q_{2}}, 
\end{equation*}
then
\begin{equation}\label{imp-ident}
\int_{\Omega} (q_{1} - q_{2}) v_{1} v_{2} = 0
\end{equation}
for any $v_{i} \in H^{1}(\Omega)$  ($i=1, \, 2$) a solution of the equation
\begin{equation*}
\Delta v_{i} - q_{i} v_{i} = 0 \mbox{ in } \Omega.
\end{equation*}
The crucial idea of Sylvester and Uhlmann in \cite{SylvesterUhlmann87} is to the find a (large) class of solutions 
of the equation
\begin{equation*}
\Delta v - q v = 0  \mbox{ in } \mR^d
\end{equation*} 
of the form 
\begin{equation*}
v = (1 + w) e^{x \cdot \xi / 2} \mbox{ in } \mR^{d}, 
\end{equation*}
where $\xi \in \mC^{d}$ with $\xi \cdot \xi  = 0$ and $|\xi|$ is large. Since $\xi \cdot \xi = 0$, it follows that 
\begin{equation}\label{iteration}
\Delta w + \xi \cdot \nabla w  - q w = q \mbox{ in } \mR^{d}; 
\end{equation}
here one extends $q$ appropriately on $\mR^{d}$ and denotes the extension also by $q$.   
Their key observation is
\begin{equation}\label{SU87-0}
\lim_{|\xi| \to 0}\| w\|_{H^{1}(B_{r})} = 0 \mbox{ for } r > 0,  
\end{equation}
which is a consequence of the following fundamental estimate established in \cite{SylvesterUhlmann87}: 
\begin{equation}\label{SU87}
\| W \|_{H^{1}(B_{r})} \le \frac{C_{r}}{|\xi|} \| f\|_{H^{1}} \quad \forall \; r > 0, 
\end{equation}
if $f$ has compact support, where $W$ is the solution to the equation
\begin{equation}\label{eq-main}
\Delta W + \xi \cdot \nabla W = f \mbox{ in } \mR^{d}.  
\end{equation}
By appropriate choices of $\xi_{1}$ and $\xi_{2}$ for the associated 
$v_{1}$ and $v_{2}$ 
with $\xi_{1} + \xi_{2} = 2k$, a constant vector in $\mC^{d}$, then
using \eqref{imp-ident} and \eqref{SU87-0}, they  show that 
\begin{equation*}
\int_{\Omega} (q_{1} - q_{2}) e^{k \cdot x} = 0 \mbox{ for all } k \in \mC^{d}. 
\end{equation*}
This in turn  implies 
\begin{equation*}
q_{1} = q_{2}.
\end{equation*} 

In \cite{BT, Brown96}, the authors improved this estimate for solutions to \eqref{eq-main} in a Besov space where $f$ has  $-1/2$ derivatives.   The proof in \cite{PaivarintaPanchenkoUhlmann03} is based on a different way of constructing CGO solutions.  

\medskip
We next discuss the approach due to Haberman and Tataru in \cite{HabermanTataru11}. The key point in \cite{HabermanTataru11} is to consider solutions to \eqref{eq-main} in $X_{\xi}^{1/2}$ with $f \in X_{\xi}^{-1/2}$, where
\begin{equation*}
\| f\|_{X_{\xi}^{s}}: = \big\| \big| |k |^{2} + k \cdot \xi \big|^{s}\hat f(k) \big\|_{L^{2}} \quad \mbox{ for } s \in \mR.  
\end{equation*}
These special function spaces have roots from the work of  Bourgain in \cite{B}. Their key estimates involves various quantities related to $L^{2}$-norm of a function by its norm in $X_{\xi}^{s}$ with $s =-1/2$ or $1/2$. This is given in   \cite[Lemma 2.2]{HabermanTataru11}.  Another ingredient in their proof is an averaging estimate for solutions to \eqref{eq-main}, \cite[Lemma 3.1]{HabermanTataru11}.  

\medskip
The work of Kenig and Jerison in \cite{Chanillo} is in the spirit of  \cite{SylvesterUhlmann87} but  uses  a  generalized Sobolev inequality, due to Kenig et al.  in \cite{KRS}.  This Sobolev inequality
 for $W$, a solution to \eqref{eq-main}, is of the form
\begin{equation}\label{Sobolev}
\| W\|_{L^{p}} \le C \| f\|_{L^{p'}}, 
\end{equation}
if $1 < p < + \infty$ and $1 < p': = pd/ (d+2) < + \infty$.  In \cite{Chanillo} the requirement $p> d/2$ is used to showed that 
\begin{equation*}
\| W\|_{L^{q}} \le C |\xi|^{-\alpha} \| V\|_{L^{q}},  
\end{equation*}
where $\alpha = 2 - d/ p$ and $(q-2)/q = 1/p$, and $W$ is the solution to equation \eqref{eq-key} below. This estimate was used in their iteration process to obtain solutions to \eqref{iteration}. 

\medskip
The construction of CGO solutions by  Ferreira et al. in \cite{FerreiraKenigSalo} is quite different and based on a limiting Carleman's estimate originating in the work of \cite{FerreiraKenigSalo1}.  

\medskip 
The goal of the paper is to introduce an approach, which is based on the construction of CGO solutions in the spirit of Sylvester and Uhlmann and involves only standard Sobolev spaces, to prove the following results: 

\begin{enumerate}
\item[i)] $\Lambda_{q}$ uniquely determines $q$ if $q = \dive g_{1} + g_{2}$ where 
$\inf_{\phi \in [C(\bar \Omega)]^d}\|g_1- \phi  \|_{L^\infty}$ is small, $g_{1} \in L^{\infty}(\Omega) \cap C^0(\overline{\Omega}_\delta)$ for some $\delta >0$,  $\hat g_{1} \in L^{p}$ for some $p<2$, and $g_{2} \in L^{d}$.  Here $\overline{\Omega}_\delta = \overline{ \{ \dist(x,\partial \Omega) < \delta \} \cap \Omega}$ (Theorem~\ref{thm1}).

\item[ii)] $\Lambda_{q}$ uniquely determines $q$  for $q \in L^{d/2}$ (Theorem~\ref{thm2}). 

\item[iii)]  $\Lambda_{q}$ uniquely determines $q$ if $q = \dive g_1 + g_2$ where $g_1 \in W^{t, 3/(t+1)}(\Omega)$ for some $t>1/2$ and $g_2 \in L^{3/2}(\Omega)$ in three dimensions (Theorem~\ref{thm3}).  
\end{enumerate}
 \medskip

To this end, we extend results of Sylvester and Uhlmann in \cite{SylvesterUhlmann87} on the stability of solutions to \eqref{eq-main} for one negative order. The proof is different from the one in \cite{SylvesterUhlmann87} and  quite simple. The same approach also  implies similar results as in  \cite{SylvesterUhlmann87}. The proof of $i)$ is mainly based on a new observation on the stability of the following equation 
(see Lemma~\ref{lem1}) 
\begin{equation}\label{eq-key}
\Delta W + \xi \cdot \nabla W = q V \mbox{ in } \mR^{d},  
\end{equation}
which is the key for the iteration process to obtain a solution to \eqref{iteration}, and an averaging argument for initial data  (see Lemma~\ref{lem-average1}) in the same spirit of \cite{HabermanTataru11}.  The (new) proof of $ii)$ in this paper is (only) based on a combination of the generalized Sobolev inequality and the standard approach used in \cite{SylvesterUhlmann87} (see Proposition~\ref{pro1-Nachmann}); however, the iteration process used to obtain solutions to \eqref{iteration} is quite tricky. The proof of $iii)$ is based on an averaging argument on {\bf both initial data and the kernel} (see Lemmas~\ref{lem-new} and \ref{lem-average-11}).

\medskip 
Statement i) is slightly different from what one can derive directly from the results of Haberman and Tataru. Statement ii) is Lavine and Nachman's result. 
Statment iii) {\bf implies the results of Lavine and Nachmann in three dimensions and yields uniqueness for a larger class of conductivities}. 
As a consequence, we  give a new proof for Haberman and Tataru's result under a mild additional assumption (Corollary~\ref{cor1}), Lavine and Nachman's result, and prove the uniqueness of Calder\'on's problem for the class of $W^{s, 3/s}$ (for some $s> 3/2$) conductivities in three dimensions (Corollary~\ref{cor3}); this last result is new as far as we know. 

\medskip
Let us describe the ideas of the proof of each conclusion in more detail. 
Without loss of generality one may assume that $\supp q \subset B_1$. Here and in what follows $B_{r}(a)$ denotes the ball centered at $a$ of radius $r, $ and $B_{r}$ denotes $B_{r}(0)$. 
Concerning i), our new key estimate  for solutions to \eqref{eq-key} is 
\begin{equation*}
\begin{split}
& 
\| \nabla W\|_{L^{2}(B_{r})} + |\xi | \cdot \| W\|_{L^{2}(B_{r})} \\
& \qquad \qquad \le C_r \big(\| g_{1}\|_{L^{\infty}}  + \| g_{2}\|_{L^{d}} \big) \Big(\| \nabla V\|_{L^{2}(B_{1})} + |\xi | \cdot \| V\|_{L^{2}(B_{1})} \Big), 
\end{split}
\end{equation*}
if $q = \dive g_{1} + g_{2}$ and $\supp g_1, \; \supp g_2 \subset \subset B_1$, see Lemma~\ref{lem1}. The proof of this inequality is based on an estimate for solutions to \eqref{eq-main} in which $f \in H^{-1}$ in the spirit \eqref{SU87} and is  presented  in Lemma~\ref{fund-lemma}. The proof of Lemma~\ref{fund-lemma} is quite elementary and different from the proof in  \cite{SylvesterUhlmann87}. 
After this, we employ some average estimates, as in \cite{HabermanTataru11}.   We remark that 
we will need $g_1 \in C^0(\overline{\Omega}_\delta)$ to ensure the existence of a trace when turning 
the elliptic PDE \eqref{eq1} into the integral \eqref{imp-ident}. 
Concerning ii), we first split $q$ into $f + g$ where $f$ is smooth and $\| g\|_{L^{d/2}}$ is small. Using the generalize Sobolev inequality \eqref{Sobolev} and the standard estimates for CGO solutions \eqref{SU87}, we are able to  reach  
\begin{equation*}
\lim_{|\xi| \to \infty} \|w \|_{H^1(B_r)}= 0
\end{equation*} 
where $w$ is the solution to \eqref{iteration}. The iteration process to obtain the existence of $w$ and the estimate of $w$ mentioned above are rather tricky in this case.  
Concerning iii), our key ingredient are 1) the following estimate for solutions to \eqref{eq-key} 
\begin{equation*}
\| W\|_{H^{1}(B_{r})} \le E(q, \xi) \| V\|_{H^{1}}
\end{equation*}
for some $E(q,\xi)$ (see Lemma~\ref{lem-new}),
and  2) the observation that, roughly speaking,  if $q \in H^{-1/2}$ with compact support then $E(q, \xi) \to 0$ as $\xi \to \infty$ for a large set of $\xi$'s (see Proposition~\ref{pro1-new}). 
At this point we both average as in \cite{HabermanTataru11} and also average $E(q, \xi)$; the estimate for solutions of \eqref{eq-key} depends on the direction of $\xi$ and $q$. 

\medskip
We state these results explicitly. Concerning i), using the construction of  CGO solutions in the spirit of Sylvester and Uhlmann in standard Sobolev spaces and some new observations (Lemma~\ref{lem1}, see also \ref{SU2}), we can reach 

\begin{theorem}\label{thm1} Let $d \ge 3$, $\Omega$ be a smooth bounded subset of $\mR^{d}$. Let $g_{1}, h_{1} \in L^{\infty}(\Omega) \cap { C^0(\overline{\Omega}_\delta)}$ for any $\delta > 0$, $g_{2}, h_{2} \in L^{d}(\Omega)$ be such that 
\begin{equation}\label{Lp}
\| {\mathcal F} (1_{\Omega} g_{1})\|_{L^{p}} + \|{\mathcal F} (1_{\Omega }h_{1}) \|_{L^{p}} < \infty
\footnote{Here $1_{\Omega}$ denotes the characteristic function of $\Omega$ and ${\mathcal F}$ denotes the Fourier transform.  This technical condition 
arises from our averaging estimates in Lemma~\ref{lem-average1} (see Remark~\ref{rem-Tech}). } \mbox{ for some }  1 < p<2.
\end{equation}  
Set 
\begin{equation*}
q_{1} = \dive g_{1} + g_{2} \quad  \mbox{ and } \quad q_{2} = \dive h_{1} + h_{2}.  
\end{equation*}
Assume that 
\begin{equation*}
\Lambda_{q_{1}} = \Lambda_{q_{2}},
\end{equation*}
then  there exists a positive constant $c$ such that if 
\begin{equation}
\inf_{\phi \in C(\bar \Omega)} \| g_1 - \phi \|_{L^\infty}  + \inf_{\phi \in C(\bar \Omega)} \| h_1 - \phi \|_{L^\infty} \le c. 
\end{equation}
then 
\begin{equation*}
q_{1} = q_{2}. 
\end{equation*}
\end{theorem}

As a consequence, we obtain the following result which is slightly weaker from the one of 
Haberman and Tataru in \cite{HabermanTataru11}. 

\begin{corollary} \label{cor1}  Let $d \ge 3$, $\Omega$ be a smooth bounded subset of $\mR^{d}$,  $\gamma_{1}, \gamma_{2} \in W^{1, \infty}(\Omega) \cap C^1(\Omega_\delta)$ for some $\delta > 0$ be such that 
\begin{equation*}
1/ \lambda \le \gamma_{1}(x), \gamma_{2}(x) \le \lambda \mbox{ for } a.e. \; x \in \Omega, 
\end{equation*}
for some $\lambda > 0$ and
 \begin{equation}\label{tech-assumption}
 {\mathcal F}\big(1_{\Omega} \nabla \ln \gamma_{i} \big) \in L^{p} \mbox{ for some $1< p<2$}.
 \end{equation}
 Assume that 
\begin{equation*}
DtN_{\gamma_{1}} = DtN_{\gamma_{2}},
\end{equation*}
then there exists a positive constant $c$ such that if 
\begin{equation}\label{small-assumption}
\inf_{\phi \in [C(\bar \Omega)]^d}\|\nabla \ln \gamma_{1} - \varphi \|_{L^{\infty}} + \inf_{\phi \in [C(\bar \Omega)]^d} \|\nabla \ln \gamma_{2} - \varphi \|_{L^{\infty}}  \le c, 
\end{equation} 
then 
\begin{equation*}
\gamma_{1}  = \gamma_{2}. 
\end{equation*}
\end{corollary}

Assumption~\eqref{tech-assumption} is a mild condition since it holds  holds for $p=2$ since $g_{1}, h_{1} \in L^{\infty}(\Omega)$. Assumption~\eqref{tech-assumption} is not required in \cite{HabermanTataru11}. The requirement that $\gamma_1, \gamma_2 \in C^1(\Omega_\delta)$  does not appear in \cite{HabermanTataru11}. Statement \eqref{small-assumption} is stronger than their results; however, their method can derive \eqref{small-assumption} as well.

\medskip
Concerning ii), we give a new proof of 

\begin{theorem} \label{thm2}
Let {$d \ge 3$}, $\Omega$ be a smooth bounded subset of $\mR^{d}$. Let $q_{1}, q_{2} \in L^{d/2}(\Omega)$.
Assume that 
\begin{equation*}
\Lambda_{q_{1}} = \Lambda_{q_{2}},
\end{equation*}
then 
\begin{equation*}
q_{1} = q_{2}. 
\end{equation*}
\end{theorem}

As a consequence of Theorem~\ref{thm2}, one obtains

\begin{corollary}\label{cor2} Let {$d \ge 3$}, $\Omega$ be a smooth bounded subset of $\mR^{d}$,  $\gamma_{1}, \gamma_{2} \in W^{2, d/2}(\Omega)$ be such that 
\begin{equation*}
1/ \lambda \le \gamma_{1}(x), \gamma_{2}(x) \le \lambda \mbox{ for } a.e. \; x \in \Omega, 
\end{equation*}
for some $\lambda > 0$. Assume that 
\begin{equation*}
DtN_{\gamma_{1}} = DtN_{\gamma_{2}},
\end{equation*}
then 
\begin{equation*}
\gamma_{1}  = \gamma_{2}. 
\end{equation*}
\end{corollary}

Concerning iii), we obtain the following new result 

\begin{theorem}\label{thm3} Let  $\Omega$ be a smooth bounded subset of $\mR^{3}$,  {$g_{1}, h_{1} \in W^{t, 3/(t+1)}(\Omega) $ for some  $t > 1/2$, $g_{2}, h_{2} \in L^{3/2}(\Omega)$. } Set 
\begin{equation*}
q_{1} = \dive g_{1} + g_{2} \mbox{ and } q_{2} = \dive h_{1} + h_{2}.  
\end{equation*}
Assume that 
\begin{equation*}
\Lambda_{q_{1}} = \Lambda_{q_{2}},
\end{equation*}
then 
\begin{equation*}
q_{1}  = q_{2}. 
\end{equation*}
\end{theorem}

Here is a consequence of Theorem~\ref{thm3}.

\begin{corollary}\label{cor3} Let  $\Omega$ be a smooth bounded subset of $\mR^{3}$, { $\gamma_{1}, \gamma_{2} \in W^{s, 3/s}(\Omega)$ for some $s>3/2$} be such that 
\begin{equation*}
1/ \lambda \le \gamma_{1}(x), \gamma_{2}(x) \le \lambda \mbox{ for } a.e. \; x \in \Omega, 
\end{equation*}
for some $\lambda > 0$. Assume that 
\begin{equation*}
DtN_{\gamma_{1}} = DtN_{\gamma_{2}},
\end{equation*}
then 
\begin{equation*}
\gamma_{1}  = \gamma_{2}. 
\end{equation*}
\end{corollary}

%

\begin{remark}\label{rem1} In Theorems~\ref{thm1}, \ref{thm2},  and \ref{thm3}, $0$ is assumed not a Dirichlet eigenvalue for the potential problems. Then the fact that $\Lambda_{q_{1}} = \Lambda_{q_{2}}$ implies $q_{1} = q_{2}$. In fact this assumption can be weaken as follows. Assume that 
\begin{equation*}
\frac{\partial v_{1}}{\partial \eta} = \frac{\partial v_{2}}{\partial \eta},   
\end{equation*}
for any $v_{1}, v_{2} \in H^{1}(\Omega)$ such that 
\begin{equation*}
\Delta v_{i} - q_{i} v_{i} = 0 \mbox{ in } \Omega \mbox{ for } i=1, 2, \mbox{ and } v_{1} = v_{2} \mbox{ on } \partial \Omega. 
\end{equation*}
Then $q_{1} = q_{2}$ under the same conditions on $q_{i}$, $i=1,2$. In fact, we prove Theorems~\ref{thm1}, \ref{thm2}, and \ref{thm3} under this weaker assumption.  
\end{remark}

The paper is organized as follows.  In Section~\ref{sec-pre}, we establish new estimates for CGO solutions in the spirit of Sylvester and Uhlmann.  In Section~\ref{sec-Tataru} we establish Theorem~\ref{thm1} and  Corollary~\ref{cor1}. This is established by generating CGO solutions via a direct iteration method and averaging methods.  We then turn to the proof of Theorem~\ref{thm2} and Corollary~\ref{cor2} in Section~\ref{lavinenachman}.  Section~\ref{newclass} handles the proof of Theorem~\ref{thm3} and Corollary~\ref{cor3}. Finally, in Appendix~\ref{apA} we provide a few results on averaging of the kernel $K_\xi(x)$ to \eqref{eq-key} that are used crucially in our CGO arguments, and in Appendix~\ref{apB} we establish that $\gamma_1 = \gamma_2$ on $\partial \Omega$ if $DtN_{\gamma_1} = DtN_{\gamma_2}$ and $\gamma_1, \gamma_2$ belong only to $W^{1, 1}(\partial \Omega)$. 
We recall again that all above results are only obtained via the construction of CGO solutions in standard Sobolev spaces and averaging arguments.


\section{New estimates for CGO solutions in the spirit of \\ Sylvester and Uhlmann}\label{sec-pre}

In this section, we recall and extend the fundamental estimates due to Sylvester and Uhlmann in \cite{SylvesterUhlmann87} concerning solutions of the equation 
\begin{equation}\label{fund-equation}
\Delta w + \xi \cdot \nabla w = f 
\end{equation}
where $\xi \in C^{d}$ and $\xi \cdot \xi = 0$. 

Given $\xi \in \mC^{d}$ with $|\xi|> 2$ and $\xi\cdot \xi = 0$, define
\begin{equation*}
\widehat{K_{\xi}}(k) = {1\over  -|k|^2 + i \xi \cdot k} \quad \mbox{ for } k \in \mR^{d}. 
\end{equation*}
Then for $f \in H^{-1}(\mR^{d})$ with compact support, $K_{\xi} * f$ is a solution to the equation 
\begin{equation*}\label{toto}
\Delta w + \xi \cdot \nabla w  = f \mbox{ in } \mR^{d},  
\end{equation*}
and 
\begin{equation*}
\widehat{ K_{\xi}*f}  = \widehat K_{\xi} \cdot \hat f \in L^{1} + L^{2}. 
\end{equation*}
We recall the following fundamental results due to Sylvester and Uhlmann in \cite{SylvesterUhlmann87}. 

\begin{proposition}[Sylvester-Uhlmann]\label{fund-lemma-SU} Let $-1 < \delta < 0$, $\xi \in \mC^{d}$ with $|\xi|> 2$ and $\xi\cdot \xi = 0$,  and let $f \in L^{2}_{\loc}(\mR^{d})$.  Then 
\begin{equation}\label{SU1-o}
\| K_{\xi} * f \|_{H^k_{\delta}} \leq {C \over |\xi|} \| f \|_{H^k_{1 + \delta}} \quad \mbox{ for } k \ge 0, 
\end{equation}
\begin{equation}\label{SU2-o}
\| K_{\xi} * f \|_{H^{k+1}_{\delta}} \leq {C} \| f \|_{H^{k}_{1 + \delta}}  \quad \mbox{ for } k \ge 0. 
\end{equation}
for some positive constant  $C$ independent of $\xi$ and $f$. 
\end{proposition}

\noindent Here 
\begin{equation*}
\| v\|_{L^{2}_{\delta}} := \| (1 + |\cdot |^{2})^{\delta} v(\cdot)\|_{L^{2}}
\end{equation*}
and 
\begin{equation*}
\| v\|_{H^{k}_{\delta}} := \sum_{|\alpha| = 0}^{k}\| (1 + |\cdot |^{2})^{\delta} D^{\alpha } v(\cdot)\|_{L^{2}}.
\end{equation*}

These estimates play an important role in their proof of the uniqueness of smooth potentials \cite{SylvesterUhlmann87} and in the proofs of the improvements in \cite{Brown96, PaivarintaPanchenkoUhlmann03, FerreiraKenigSalo}.  

\medskip
We will extend the above results to negative derivatives and to the case with two derivative difference, which are crucial for the proof of Theorem~\ref{thm1}. Our proof for negative derivatives and the  two derivative difference is rather elementary. The same proof also gives the following estimates, for $f \in L^{2}(\mR^d)$ with compact support, 
\begin{equation}\label{SU1}
\| K_{\xi} * f \|_{H^{k}(B_{r})} \leq {C_{r} \over |\xi|} \| f \|_{H^k} \quad \mbox{ for } k \ge 0, 
\end{equation}
and
\begin{equation}\label{SU2-1}
\| K_{\xi} * f \|_{H^{k+1}(B_{r})} \leq {C_{r}} \| f \|_{H^{k}}  \quad \mbox{ for } k \ge 0
\end{equation}
Here  $C$ is a positive number independent of $\xi$ and $f$. These estimates are slightly weaker than the original ones of Sylvester and Uhlmann in \eqref{SU1-o} and \eqref{SU2-o}; however, they are sufficient for establishing the uniqueness of smooth potential in \cite{SylvesterUhlmann87}. Here is the extension:

\begin{lemma}\label{fund-lemma} Let $R>0$, $\xi \in \mC^{d}$ with $|\xi|> 2$ and $\xi \cdot  \xi = 0$,  and let $f \in H^{-1}(\mR^{d})$ with $\supp f \subset B_{R}$.  Then  
\begin{equation}\label{SU2}
\| K_{\xi} * f \|_{L^{2}(B_{r}) } \le C_{r}  \| f \|_{H^{-1}} 
\end{equation}
and
\begin{equation}\label{SU3}
\| K_{\xi} * f \|_{H^{k+1}(B_{r})} \le C_{r} |\xi|  \cdot \| f \|_{H^{k-1}}, \quad \mbox{ for } k \ge 0, 
\end{equation}
for some $C_r$ which depends on $r$ and $R$ but is independent of $\xi$ and $f$. 
\end{lemma}


\medskip

\noindent{\bf Proof.} We will prove \eqref{SU2}; the proof of \eqref{SU3} as well \eqref{SU1} and \eqref{SU2-1} follow similarly. Set
\begin{equation*}
\Gamma_{\xi} := \{k \in \mR^{d} ; \; - |k|^{2}  + i \xi \cdot k =0 \}. 
\end{equation*}
It is clear that 
\begin{equation}\label{pro-K1}
|\hat K_{\xi}(k)| \le \frac{C}{|\xi| \dist (k, \Gamma_{\xi})} \mbox{ if } |k| \le 2 |\xi|, \mbox{ and } |\hat K_{\xi}(k)|  \le \frac{C}{|k|^{2}} \mbox{ if } |k| \ge 2|\xi|,
\end{equation}
In this proof, $C$ denotes a positive constant independent of $\xi$ and $f$. 
Define $K_{1, \xi}$ and $K_{2, \xi}$ as follows
\begin{equation}\label{def-K1}
\hat K_{1, \xi} (k) = 
\left\{\begin{array}{cl}
\hat K_{\xi}(k) &  \mbox{ if } \dist(k, \Gamma_\xi) \ge 1,  \\[6pt]
0 &  \mbox{ otherwise}, 
\end{array} \right. 
\end{equation}
and 
\begin{equation}\label{def-K2}
\hat K_{2, \xi} (k) = \hat K_{\xi}(k) - \hat K_{1, \xi}(k),
\end{equation}
and so
\begin{equation}\label{b1}
\| K_{\xi} * f \|_{L^{2}(B_{r})} \le \| K_{1, \xi} * f \|_{L^{2}(B_{r})}  + \| K_{2, \xi} * f \|_{L^{2}(B_{r})}. 
\end{equation}
Using Plancherel's theorem, we derive from \eqref{pro-K1} and \eqref{def-K1} that
\begin{equation}\label{b2}
\| K_{1, \xi} * f \|_{L^{2}(\mR^{d})} \le C \| f\|_{H^{-1}}. 
\end{equation}
Fix 
\begin{equation}\label{choice-phi}
\varphi \in C^{\infty}_{0}(\mR^{d}) \mbox{ with } \varphi = 1 \mbox{ in } B_{R + r}.
\end{equation} 
Since $\supp f \subset B_{R}$, it follows that $f =  \varphi f$; hence
\begin{equation*}
\hat f = \hat \varphi * \hat f. 
\end{equation*}
Define 
\begin{equation}\label{def-f1}
\widetilde f(k) = \sup_{\eta \in B_{4}(k)} |\hat f(\eta)|
\end{equation}
and 
\begin{equation}\label{def-phi1}
\widetilde \varphi (k) = \sup_{\eta \in B_{4}(k)} |\hat \varphi (\eta) | .
\end{equation}
Since 
\begin{equation*}
 |\hat f|* |\hat \varphi| (\eta) = \int_{\mR^{d}} |\hat f (\zeta)| |\hat \varphi(\eta -\zeta)| \, d \zeta, 
\end{equation*} 
it follows from the definition of $\widetilde f$ \eqref{def-f1} and  $\widetilde \varphi$ \eqref{def-phi1} that 
\begin{equation}\label{bound-f}
\widetilde f \le |\hat f|* \widetilde \varphi. 
\end{equation}
From the choice of $\varphi$ \eqref{choice-phi}, we have 
\begin{equation}\label{to0}
\begin{split}
\| K_{2, \xi} * f \|_{L^{2}(B_{r})}^{2} & \le \|\varphi (K_{2, \xi} * f) \|_{L^{2}(\mR^{d})}^{2} \\
&  \le  \int_{\mR^{d}} \Big| \int_{\dist(\eta, \Gamma_{\xi}) \le 1} |\hat \varphi(k- \eta)| \cdot |\hat K_{\xi} (\eta)| \cdot |\hat f(\eta)| \, d \eta \Big|^{2} \, dk. 
\end{split}
\end{equation}
Using the fact that 
\begin{equation}\label{pro-K}
\int_{|x| \le 1} \frac{1}{|x_{1}| + |x_{2}|} \, dx_1 \, d x_2 < + \infty, 
\end{equation}
we obtain
\begin{equation}\label{to1}
\int_{\dist(\eta, \Gamma_{\xi}) \le 1} |\hat \varphi(k- \eta)| \cdot  |\hat K_{\xi} (\eta)| \cdot |\hat f(\eta)| \, d \eta \le  \frac{C}{|\xi|} \int_{\dist(\eta, \Gamma_{\xi}) \le 1} \widetilde \varphi(k- \eta) \widetilde f(\eta) \, d \eta. 
\end{equation}
In fact, for $|\xi| > 2$, there exists $0 < r \le 1$ (independent of $\xi$) such that for $\eta$ with $\dist(\eta, \Gamma_\xi) \le r$, there exists an unique pair $(\eta_1, \eta_2) \in \mR^{d} \times \mR^d$ such that 
$\eta_1 \in \Gamma_\xi$, $\eta_2 \perp T_{\Gamma_\xi}(\eta_1)$, the tangent plane of $\Gamma_\xi$ at $\eta_1$, such that $|\eta_2| \le r$ and $\eta_1 + \eta_2 = \eta$. Then 
\begin{multline*}
\int_{\dist(\eta, \Gamma_{\xi}) \le r} |\hat \varphi(k- \eta)| \cdot  |\hat K_{\xi} (\eta)| \cdot |\hat f(\eta)| \, d \eta \\[6pt]
\le C \int_{\eta_1 \in \Gamma_\xi} \int_{|\eta_2| \le r; \eta_2 \perp T_{\Gamma_\xi}(\eta_1)} |\hat \varphi(k- \eta_1 - \eta_2)| \cdot  |\hat K_{\xi} (\eta_1 + \eta_2)| \cdot |\hat f(\eta_1 + \eta_2)| \, d \eta_2 \, d \eta_1. 
\end{multline*}
Since 
\begin{multline}
\int_{\eta_1 \in \Gamma_\xi} \int_{|\eta_2| \le r; \eta_2 \perp T_{\Gamma_\xi}(\eta_1)} |\hat \varphi(k- \eta_1 - \eta_2)| \cdot  |\hat K_{\xi} (\eta_1 + \eta_2)| \cdot |\hat f(\eta_1 + \eta_2)| \, d \eta_2 \, d \eta_1 \\[6pt]
\le 
\int_{\eta_1 \in \Gamma_\xi} \sup_{|\eta_2| \le r} |\hat \varphi(k - \eta_1 - \eta_2) | \sup_{|\eta_2| \le r} |\hat f(\eta_1 + \eta_2)|  \int_{|\eta_2| \le r; \eta_2 \perp T_{\Gamma_\xi}(\eta_1)}  |\hat K_{\xi} (\eta_1 + \eta_2)|  \, d \eta_2 \, d \eta_1. 
\end{multline}
and, by \eqref{pro-K}, 
\begin{equation*}
 \int_{|\eta_2| \le r; \eta_2 \perp T_{\Gamma_\xi}(\eta_1)}  |\hat K_{\xi} (\eta_1 + \eta_2)|  \, d \eta_2\le \frac{C}{|\xi|},  
\end{equation*}
it follows that 
\begin{equation}\label{ha1}
\int_{\dist(\eta, \Gamma_{\xi}) \le r} |\hat \varphi(k- \eta)| \cdot  |\hat K_{\xi} (\eta)| \cdot |\hat f(\eta)| \, d \eta \le \frac{C}{|\xi|} \int_{\eta_1 \in \Gamma_\xi} \sup_{|\eta_2| \le r} |\hat \varphi(k - \eta_1 - \eta_2) | \sup_{|\eta_2| \le r} |\hat f(\eta_1 + \eta_2)| \, d \eta_1. 
\end{equation}
On the other hand, by the definition of $\tilde f$ and $\tilde \varphi$, 
\begin{equation}\label{ha2}
\int_{\eta_1 \in \Gamma_\xi} \sup_{|\eta_2| \le r} |\hat \varphi(k - \eta_1 - \eta_2) | \sup_{|\eta_2| \le r} |\hat f(\eta_1 + \eta_2)| \, d \eta_1 \le C \int_{\dist(\eta, \Gamma_{\xi}) \le 1} \widetilde \varphi(k- \eta) \widetilde f(\eta) \, d \eta. 
\end{equation}
A combination of \eqref{ha1} and \eqref{ha2} yields \eqref{to1}. 

\medskip

Applying H\"older's inequality, we derive from \eqref{to0} and \eqref{to1} that  
\begin{equation}\label{b3-1}
\| K_{2, \xi} * f \|_{L^{2}(B_{r})}^{2} \le  \frac{C}{|\xi|^{2}} \int_{\dist(\eta, \Gamma_{\xi}) \le 1} |\widetilde f (\eta)|^{2} \,  d \eta. 
\end{equation}
We now estimate the RHS of \eqref{b3-1}. 
Applying H\"older's inequality again, from  \eqref{bound-f} and the fact that $\tilde \varphi \in L^1$, we have
\begin{equation}\label{sss}
 \int_{\dist(\eta, \Gamma_{\xi}) \le 1} |\widetilde f (\eta)|^{2} \, d \eta \le C \int_{\dist(\eta, \Gamma_{\xi}) \le 1} \int_{\mR^{d}} \widetilde \varphi(\eta - k) |\hat f(k)|^{2} \, dk \, d \eta.  
\end{equation}
Using Fubini's theorem, we derive from \eqref{sss} that 
\begin{equation}\label{b3-2-1}
 \int_{\dist(\eta, \Gamma_{\xi}) \le 1} |\widetilde f (\eta)|^{2} \, d \eta \le C \int_{\mR^{d}} |\hat f(k)|^{2}  \int_{\dist(\eta, \Gamma_{\xi}) \le 1} \widetilde \varphi(\eta - k) \, d \eta  \, dk. 
\end{equation}
Since $\tilde \varphi \in {\mathcal S}$, the Schwartz class, it follows that 
\begin{equation}\label{b3-2-2}
 \int_{\mR^{d}} |\hat f(k)|^{2}  \int_{\dist(\eta, \Gamma_{\xi}) \le 1} \widetilde \varphi(\eta - k) \, d \eta  \, dk
 \le C \Big( \int_{|k| \le 2 |\xi|} |\hat f(k)|^{2} \, dk  + 
 \int_{|k| > 2 |\xi|} \frac{|\hat f(k)|^{2}}{|k|^{2}} \, dk  \Big). 
 \end{equation}
From \eqref{b3-2-1} and \eqref{b3-2-2}, we obtain
\begin{equation}\label{b3-2}
\frac{1}{|\xi|^{2}} \int_{\dist(\eta, \Gamma_{\xi}) \le 1} |\widetilde f (\eta)|^{2} \, d \eta \le C \|f\|_{H^{-1}}^{2}.
\end{equation}
A combination of  \eqref{b3-1} and \eqref{b3-2} yields, 
\begin{equation}\label{b3}
\| K_{2, \xi} * f \|_{L^{2}(B_{r})}^{2} \le  C\| f \|_{H^{-1}}^{2}. 
\end{equation}
The conclusion follows from \eqref{b1}, \eqref{b2}, and \eqref{b3}. \proofend

\section{Proof of Theorem~\ref{thm1} and Corollary~\ref{cor1}} \label{sec-Tataru}


In this section, we prove Theorem~\ref{thm1} and Corollary~\ref{cor1}. The proof of Theorem~\ref{thm1} contains two main ingredients. The first one is  a new useful inequality (Lemma~\ref{lem1})  and its variant (Lemma~\ref{lem2})
to solutions to \eqref{fund-equation} whose the proof is based on estimates presented in Section~\ref{sec-pre}.  The second one is an averaging estimate (Lemma~\ref{lem-average1}) with respect to $\xi$ for $K_{\xi}*q$ in the same spirit in   \cite{HabermanTataru11} and is presented in Appendix \ref{apA}. 

\subsection{Some useful lemmas}

The following lemma is new and interesting in itself. It plays an important role in our analysis. Its proof is quite elementary, and can be seen as the replacement of  \cite[Lemma 2.3]{HabermanTataru11}.

\begin{lemma}\label{lem1} Let $d \ge 3$, $\xi \in \mC^{d}$ $(|\xi|> 2)$ with $\xi \cdot \xi = 0$, $g_{1} \in [L^{\infty}(\mR^{d})]^{d}$, $g_{2} \in L^{d}(\mR^{d})$ and $V \in H^{1}(\mR^{d})$ be such that $\supp g_{1}, \supp g_{2} \subset B_{1}$. 
Set 
\begin{equation*}
q = \dive g_{1} + g_{2}
\end{equation*}
and define
\begin{equation*}
W = K_{\xi} * (qV).
\end{equation*}
We have
\begin{equation}\label{claim1}
\begin{split}
& \| \nabla W \|_{L^{2}(B_{r})} + |\xi| \cdot \| W\|_{L^{2}(B_{r})} \\
& \qquad \qquad\qquad \le C_{r} \big(\| g_{1}\|_{L^{\infty}} + \| g_{2}\|_{L^{d}} \big) \big( \| \nabla V \|_{L^{2}} + |\xi| \cdot \| V\|_{L^{2}}\big), 
\end{split}
\end{equation}
for some positive constant $C_{r}$ independent of $\xi$, $g_{1}$, $g_{2}$, and $v$. 
\end{lemma}


\noindent{\bf Proof.} We have
\begin{equation}\label{decomposition}
q V =  \dive (V g_{1} ) - g_{1} \cdot \nabla V + g_{2} V \mbox{ in } \mR^{d}. 
\end{equation}
Applying \eqref{SU2-1} with $k=1$ and \eqref{SU3} with $k=0$, we have
\begin{equation*}
\| \nabla W\|_{L^{2}(B_{r})} \le  C_{r} \Big( | \xi| \cdot \| \dive (V g_{1}) \|_{H^{-1}} 
+  \| g_{1} \cdot \nabla V \|_{L^{2}}  + \|g_{2} V\|_{L^{2}} \Big),
\end{equation*}
which implies
\begin{equation*}
\| \nabla W\|_{L^{2}(B_{r})} \le  C_{r} 
\Big(|\xi| \cdot  \|V g_{1}\|_{L^{2}} + \| g_{1} \cdot \nabla V \|_{L^{2}} + \|g_{2} V\|_{L^{2}} \Big). 
\end{equation*}
It follows that 
\begin{equation}\label{est1}
\| \nabla W\|_{L^{2}(B_{r})} \le  C_{r} \big(\| g_{1}\|_{L^{\infty}} + \| g_{2}\|_{L^{d}}\big) \big( |\xi|  \cdot \|V \|_{L^{2}}   + \| \nabla V\|_{L^{2}} \big). 
\end{equation}
Similarly, using \eqref{decomposition} and applying \eqref{SU1} with $k=0$, and \eqref{SU2}, we obtain
\begin{equation*}
|\xi| \cdot \| W\|_{L^{2}(B_{r})}  \le C_{r} \Big( |\xi| \cdot \| \dive (g_{1} V) \|_{H^{-1}} 
+  \| g_{1} \nabla V \|_{L^{2}}  + \|g_{2} V\|_{L^{2}} \Big),  
\end{equation*}
which implies
\begin{equation*}
|\xi| \cdot \| W\|_{L^{2}(B_{r})} \le C_{r} \Big(|\xi| \cdot  \|g_{1} V \|_{L^{2}} + \| g_{1} \nabla V \|_{L^{2}} + \|g_{2} V\|_{L^{2}}\Big).
 \end{equation*}
It follows that 
\begin{align}\label{est2}
|\xi| \cdot \| W\|_{L^{2}(B_{r})} \le C_{r} \big(\| g_{1}\|_{L^{\infty}} + \| g_{2}\|_{L^{d}}\big)  \big( |\xi| \cdot \|V \|_{L^{2}}   + \| \nabla V\|_{L^{2}} \big). 
\end{align}
A combination of \eqref{est1} and \eqref{est2} yields \eqref{claim1}.  \proofend

\medskip
When $g_1$ and $g_2$ are smooth, we can improve the conclusion in Lemma~\ref{lem1} as follows. 

\begin{lemma}\label{lem2} Let $d \ge 3$, $\xi \in \mC^{d}$ $(|\xi|> 2)$ with $\xi \cdot \xi = 0$, $g_{1} \in [C^{2}(\mR^{d})]^{d}$, $g_{2} \in C^{1}(\mR^{d})$ with $\supp g_{1}, \supp g_{2} \subset B_{1}$, and let $V \in H^{1}(\mR^{d})$. Set 
\begin{equation*}
q = \dive g_{1} + g_{2} 
\end{equation*}
and define 
\begin{equation*}
W = K_{\xi} * (qV). 
\end{equation*}
We have, for $r>0$, 
\begin{equation}\label{claim2}
\begin{split}
& \| \nabla W\|_{L^{2}(B_{r})} + |\xi| \cdot \| W\|_{L^{2}(B_{r})} \\
& \qquad \qquad \qquad \le \frac{C_{r}}{|\xi|} \Big(\| g_{1}\|_{C^{2}} + \| g_{2}\|_{C^{1}} \Big)  \Big( \| \nabla V \|_{L^{2}} + |\xi|  \cdot \| V\|_{L^{2}}  \Big). 
\end{split}
\end{equation}
Here $C_{r}$ is a positive constant depending only on $r$ and $d$. 
\end{lemma}

\noindent{\bf Proof.} Applying \eqref{SU1} with $k = 1$, we have
\begin{align}
 \| \nabla W\|_{L^{2}(B_{r})}&   \le   \frac{C_r}{|\xi|} \Big( \| V \dive g_{1}  \|_{H^{1}} 
+ \|V g_{2} \|_{H^{1}} \Big)  \nonumber \\
&  \le \frac{C_r}{|\xi|}  \| V \|_{H^{1}} \big( \| g_{1}\|_{C^{2}} + \|g_{2} \|_{C^{1}} \big).\label{claim2-1}
\end{align}
Similarly, 
\begin{align}
 \| W\|_{L^{2}(B_{r})} &  \le  \frac{C_r}{|\xi|} \Big( \| V \dive g_{1}  \|_{L^{2}} 
+ \|V g_{2} \|_{L^{2}} \Big)  \nonumber \\
&   \le \frac{C_r}{|\xi|}  \| V \|_{L^{2}} \big( \| g_{1}\|_{C^{1}} + \|g_{2} \|_{C^{0}} \big).\label{claim2-2}
\end{align} 
A combination of \eqref{claim2-1} and \eqref{claim2-2} yields \eqref{claim2}. \proofend

\subsection{Construction of CGO solutions}

We begin this section with an estimate for the solution of the equation 
\begin{equation*}
\Delta w + \xi \cdot \nabla w - q w = q \mbox{ in } \mR^{d}. 
\end{equation*}

\begin{proposition}\label{pro1} Let $\xi \in \mC^{d}$ $(|\xi|> 2)$ with $\xi \cdot \xi = 0$, $g_{1} \in [L^{\infty}(\mR^{d})]^{d}$, $g_{2} \in L^{d}(\mR^{d})$ with $\supp g_{1}, \supp g_{2} \subset B_{1}$. Set $q = \dive g_{1} + g_{2}$. Then there exists a  positive constant $c$ such that  if 
\begin{equation*} 
\inf_{\phi \in [C(\mR^d)]^d, \supp \phi \subset B_1 } \| g_{1} - \phi\|_{L^{\infty}}  \le c, 
\end{equation*}
then there exists 
$w \in H^{1}_{\loc}(\mR^{d})$  such that 
\begin{equation*}
w = K_{\xi} * (q  + q w)
\end{equation*}
and
\begin{equation}\label{state2-pro1}
\begin{split}
& \| \nabla (w - K_{\xi} * q)\|_{L^{2}(B_{r})} + |\xi| \cdot \| w - K_{\xi} * q\|_{L^{2}(B_{r})} \\
& \qquad \qquad  \le C_{r} \Big(\| \nabla K_{\xi} * q\|_{L^{2}(B_{2})} + |\xi| \cdot \| K_{\xi} * q\|_{L^{2}(B_{2})}\Big) \quad \forall \;  r > 0,  
\end{split}
\end{equation}
for $|\xi|$ large enough \footnote{The largeness of $|\xi|$ depends only on $g_{1}$ and $g_{2}$.}.
\end{proposition}

\noindent{\bf Proof.}   Let   $g_{i, j}$,  $1 \le i, j \le 2$,  such that 
\begin{equation*}
g_{1, 1} + g_{1, 2} = g_{1} \mbox{ and } g_{2, 1} + g_{2, 2} = g_{2}. 
\end{equation*}
\begin{equation*}
g_{1, 2}, \; g_{2, 2} \mbox{ are smooth with compact support in $B_{1}$}, 
\end{equation*}
\begin{equation*}
\|g_{1, 1}\|_{L^{\infty}} + \|g_{2, 1} \|_{L^{d}} \le 2 c,  
\end{equation*}
Set 
\begin{equation*}
q = q_{1} + q_{2},
\end{equation*}
where 
\begin{equation*}
q_{1} = \dive g_{1,1} + g_{2, 1}
\end{equation*}
and 
\begin{equation*}
q_{2} = \dive g_{1,2} + g_{2, 2}. 
\end{equation*}

Let $u_{0} = 0$  and consider the following iteration process: 
\begin{equation}\label{def-un}
w_{n}= K_{\xi} * (q + q w_{n-1}) \quad \mbox{ for } n \ge 1,
\end{equation}
which implies
\begin{equation*}
\Delta w_{n} + \xi \cdot \nabla w_{n} = q + q w_{n-1} \mbox{ in $\mR^{d}$,  for } n \ge 1.
\end{equation*}
Define 
\begin{equation*}
w_{1, n} = K_{\xi} * (q_{1} + q_{1} w_{n-1}) \quad \mbox{ and } \quad w_{2, n} = K_{\xi} * (q _{2} + q_{2} w_{n-1}). 
\end{equation*}
Then
\begin{equation*}
\Delta w_{1, n} + \xi \cdot \nabla w_{1, n} = q_{1} + q_{1} w_{n-1} \mbox{ in } \mR^{d}, 
\end{equation*}
\begin{equation*}
\Delta w_{2, n} + \xi \cdot \nabla  w_{2, n} = q _{2} + q_{2} w_{n-1} \mbox{ in } \mR^{d}, 
\end{equation*}
and 
\begin{equation}\label{ident}
w_{n} = w_{1, n} + w_{2, n} \mbox{ in } \mR^d. 
\end{equation}
Set 
\begin{equation*}
W_{n+1} = w_{n+1} - w_{n}, \quad W_{1, n+1} = w_{1, n+1} - w_{1, n}, \quad W_{2, n+1} = w_{2, n+1} - w_{2, n}. 
\end{equation*}
It follows from Lemma \ref{lem1} that 
\begin{multline}\label{id1}
\| \nabla W_{1, n+1} \|_{L^{2}(B_{r})} + |\xi| \cdot \| W_{1, n+1}  \|_{L^{2}(B_{r})} \\[6pt]
\le C_{r} \big(\|g_{1,1} \|_{L^{\infty}} + \|g_{2, 1} \|_{L^{d}}  \big)\Big( \| \nabla W_{n}\|_{L^{2}(B_{2})} + |\xi| \cdot \| W_{n}\|_{L^{2}(B_{2})} \Big),  
\end{multline}
and from Lemma \ref{lem2} that
\begin{multline}\label{id2}
\| \nabla W_{2, n+1} \|_{L^{2}(B_{r})} + |\xi| \cdot  \| W_{2, n+1} \|_{L^{2}(B_{r})} \\[6pt]
\le \frac{C_{r}}{|\xi|}  \big(\|g_{1,2} \|_{C^{2}} + \|g_{2, 2} \|_{C^{1}}  \big) \Big( \| \nabla W_{n} \|_{L^{2}(B_{2})} + |\xi| \cdot \| W_{n}\|_{L^{2}(B_{2})} \Big). 
\end{multline}
A combination of  \eqref{ident}, \eqref{id1}, and \eqref{id2} yields  
\begin{multline}\label{id3}
\| \nabla W_{n+1} \|_{L^{2}(B_{r})} + |\xi| \cdot \| W_{n+1}\|_{L^{2}(B_{r})} \\[6pt]
\le C_{r}\Big( \big(\|g_{1,1} \|_{L^{\infty}} + \|g_{2, 1} \|_{L^{d}}  \big)  
 + \frac{1}{|\xi|} \big(\|g_{1,2} \|_{C^{2}} + \|g_{2, 2} \|_{C^{1}} \big)  \Big)\Big( \| \nabla W_{n}\|_{L^{2}(B_{2})} + |\xi | \cdot  \| W_{n}\|_{L^{2}(B_{2})} \Big). 
\end{multline}
Choose $c$ such that 
\begin{equation*}
c C_{2}  = 1/ 2. 
\end{equation*}
Thus, if $|\xi|$ is large enough, then
\begin{equation*}
C_{2}\Big( \big(\|g_{1,1} \|_{L^{\infty}} + \|g_{2, 1} \|_{L^{d}}  \big)  
 + \frac{1}{|\xi|} \big(\|g_{1,2} \|_{C^{2}} + \|g_{2, 2} \|_{C^{1}} \big)  \Big) \le 3/4.
\end{equation*}
Hence, by a standard fixed point argument, it follows that 
\begin{equation*}
w_{n}  \to w \mbox{ in }  H^{1}(B_{2}).  
\end{equation*}
This implies, by  \eqref{id3}, 
\begin{equation*}
w_{n}  \to w \mbox{ in }  H^{1}(B_{r}) \quad \mbox{ for all } r > 0,    
\end{equation*}
and by \eqref{def-un}
\begin{equation*}
w= K_{\xi} * (q + q w). 
\end{equation*}
We derive from \eqref{id3} that 
 \begin{equation*}
\| \nabla  (w - w_{1})\|_{L^{2}(B_{2})} + |\xi | \cdot  \| w - w_{1} \|_{L^{2}(B_{2})} \le C \big(\| \nabla  w_{1}\|_{L^{2} (B_{2})} + |\xi |  \cdot \| w_{1} \|_{L^{2}(B_{2})}\big). 
\end{equation*}
Statement \eqref{state2-pro1} now follows from \eqref{id3}. The proof is complete. 
\proofend

\medskip 
To obtain some appropriate estimate for $u_1$ in Proposition~\ref{pro1}, we use an averaging argument in the same spirit  in \cite{HabermanTataru11}. 
More precisely, we have the following lemma whose proof is given in the appendix. 

\begin{lemma}\label{lem-average1} Let $d \ge 3$, $s> 2$, $k \in \mR^{d}$ with $|k| \ge 2$,   $1 \le p < 2$, and $R > 10$. We have 
\begin{equation}\label{est1-lem-average1}
{1\over R} \int_{R/2}^{2R} \int_{\sigma_{1} \in \mS^{d-1}} \int_{\sigma_{2} \in \mS^{d-1}_{\sigma_{1}}} |\hat K_{s\sigma_{2} - i s \sigma_{1}} (k)|^{p} \, d \sigma_{2} \, d \sigma_{1} \, ds \le C  \min \Big\{ \frac{1}{R^{p}|k|^{p}},\frac{1}{ |k|^{2p}} \Big\},  
 \end{equation}
and 
\begin{multline}\label{est2-lem-average1}
{1\over R} \int_{R/2}^{2R}  \int_{\sigma_{1} \in \mS^{d-1}} \int_{\sigma_{2} \in \mS^{d-1}_{\sigma_{1}}} \int_{\sigma_{3} \in \mS^{d-1}_{\sigma_{1}, \sigma_{2}}} \Big|\hat K_{\frac{s^2\sigma_{2}}{ \sqrt{1 + s^{2}}} + \frac{s\sigma_{3}}{ \sqrt{1 + s^{2}} }  - i s \sigma_{1} } (k) \Big|^{p}  \, d \sigma_{3}\, d \sigma_{2} \, d \sigma_{1} \, ds \\[6pt]
 \le  C \min\Big\{\frac{1}{R^{p}|k|^{p}}, \frac{1}{|k|^{2p}} \Big\},
\end{multline}
for some positive constant $C$ depending only on $d$ and $p$. 
Here  
\begin{equation}\label{S1}
\mS^{d-1}_{\sigma_{1}} := \big\{ \sigma \in \mS^{d-1}; \; \sigma \cdot \sigma_{1} = 0 \big\}
 \end{equation}
 and 
\begin{equation}\label{S12}
\mS^{d-1}_{\sigma_{1}, \sigma_{2}} := \big\{ \sigma \in \mS^{d-1}; \; \sigma \cdot \sigma_{1} = 0 \mbox{ and } \sigma \cdot \sigma_{2} = 0 \big\}. 
 \end{equation} 
\end{lemma}

\begin{remark} Let $\sigma_{1} \in \mS^{d-1}$, $\sigma_{2} \in \mS^{d-1}_{\sigma_{1}}$, and $\sigma_{3} \in \mS^{d-1}_{\sigma_{1}, \sigma_{2}}$.  Set 
\begin{equation*}
\xi_{1}  = s \sigma_{2} - i  s \sigma_{1} \quad \mbox{ and } \quad \xi_{2} =  -  \frac{s^2 \sigma_{2}}{ \sqrt{1 + s^{2}}} + \frac{s \sigma_{3}}{ \sqrt{1 + s^{2}}}  +  i s \sigma_{1},
\end{equation*}
then 
\begin{equation*}
\xi_{1} \cdot \xi_{1} =  \xi_{2} \cdot \xi_{2} = 0 \mbox{ and } \xi_{1} + \xi_{2} = \Big(s - \frac{s^2}{\sqrt{1 + s^2}} \Big) \sigma_2 + \frac{s \sigma_3}{\sqrt{1 + s^2}} \to \sigma_{3}
\end{equation*}
uniformly with respect to $\sigma_{1}$ and $\sigma_{2}$ as $s \to \infty$. 
\end{remark}

\begin{remark}\label{rem-Tech} Lemma~\ref{lem-average1} does not hold for $p=2$. The requirements \eqref{Lp} and \eqref{tech-assumption} in Theorem~\ref{thm1} and Corollary~\ref{cor1} are due to this point. 

\end{remark}

Using Proposition~\ref{pro1} and Lemma~\ref{lem-average1}, we can prove the following result: 

\begin{proposition}\label{pro2} Let  $d \ge 3$, $g_{1},  h_{1}  \in [L^{\infty}(\mR^{d})]^{d}$, $g_{2}, h_{2} \in L^{d} (\mR^{d})$ with supports in  $B_{1}$ be such that  $\hat g_{1}, \hat h_{1} \in L^{p}(\mR^{d})$ for some $1 <  p <  2$.  Assume that 
\begin{equation*}
 \inf_{\phi \in [C(\bar \Omega)]^d}\| g_{1} - \phi\|_{L^{\infty}} +\inf_{\phi \in [C(\bar \Omega)]^d} \| h_{1} - \phi \|_{L^{\infty}}  \le c, 
\end{equation*}
where $c$ is the constant in Proposition~\ref{pro1}, or $g_{1}, h_{1}$ are continuous. Set 
\begin{equation*}
q_{1} = \dive g_{1} + g_{2} \mbox{ and } q_{2} = \dive h_{1} + h_{2}. 
\end{equation*}
Then for any $0 < \eps < 1$, $n > 2$ large enough,  and $\sigma \in \mS^{d-1}$, there exist $\sigma_{1, \eps}, \, \sigma_{2, \eps}, \, \sigma_{3, \eps} \in \mS^{d-1}$, $s_{\eps} \in (n, 4n)$,  and $w_{1, \eps}, w_{2, \eps}  \in H^{1}_{\loc}(\mR^{d})$  such that 
\begin{equation}\label{orthogonality}
\sigma_{1, \eps} \cdot \sigma_{2, \eps} = \sigma_{1,\eps} \cdot \sigma_{3, \eps} = \sigma_{2, \eps} \cdot \sigma_{3, \eps} = 0, 
\end{equation}
\begin{equation}\label{dist}
|\sigma_{3, \eps} - \sigma | \le \eps, 
\end{equation}
\begin{equation*}
w_{j, \eps} = K_{\xi_{j, \eps}} * (q_{j}  + q_{j} w_{j, \eps}) \quad \mbox{ for } j =1, 2, 
\end{equation*}
and 
\begin{equation}\label{energy12}
\| \nabla w_{j, \eps}\|_{L^{2}(B_{r})} +s_{\eps} \| w_{j, \eps}\|_{L^{2}(B_{r})} \le C_{r}/ \eps^{3d} \quad \mbox{ for } j =1, 2, 
\end{equation}
for some $C_{r} > 0$ independent of $\eps$, $s$, and $\sigma$. Here
\begin{equation}\label{def-xi}
\xi_{1, \eps}  =  s_{\eps} \sigma_{2, \eps}  - i s_{\eps} \sigma_{1, \eps} \quad \mbox{ and } \quad 
\xi_{2, \eps}  =  - \frac{ s_{\eps}^2 \sigma_{2, \eps }}{ \sqrt{1 + s_{\eps}^{2}}} +  \frac{s_\eps \sigma_{3, \eps}}{ \sqrt{1 + s_{\eps}^{2}}}  + i s_{\eps} \sigma_{1, \eps}. 
\end{equation}
\end{proposition}

\noindent{\bf Proof.} Applying  Lemma~\ref{lem-average1}, we have
\begin{multline*}
\frac{1}{n}\int_{n}^{4n} \int_{\sigma_{1} \in \mS^{d-1}} \int_{\sigma_{2} \in \mS^{d-1}_{\sigma_{1}}} \Big(  |\hat K_{ s \sigma_{2} - is \sigma_{1}}(k)|^{p} \\[6pt]+   \int_{\sigma_{3} \in \mS^{d-1}_{\sigma_{1}, \sigma_{2}}} \Big|\hat K_{- \frac{ s^2\sigma_{2}}{ \sqrt{1 + s^{2}}} + \frac{s \sigma_{3}}{ \sqrt{1 + s^{2}} }  + i s \sigma_{1} } (k) \Big|^{p}  \, d \sigma_{3} \Big) \, d \sigma_{2} \, d \sigma_{1} \, ds 
 \le C \min \Big\{ \frac{1}{ n^{p} |k|^{p}}, \frac{1}{ |k|^{2p}} \Big\}. 
\end{multline*}
This implies \eqref{orthogonality} and \eqref{dist} hold for some $s_{\eps} \in (n, 4n)$ and $\sigma_{1, \eps}, \sigma_{2, \eps}, \sigma_{3, \eps} \in \mS^{d-1}$,   and 
\begin{multline}\label{energy3}
\int_{\mR^{d}} \Big( |\hat K_{\xi_{1, \eps}} (k)|^{p} |\hat q_{1}(k)|^{p} + |\hat K_{\xi_{2, \eps}} (k)|^{p} |\hat q_{2}(k)|^{p} \Big) (|k|^p + n^p) \, dk \\[6pt]
\le \frac{C}{ \eps^{3d}} \int_{\mR^{d}}  \frac{1}{ |k|^{p}} \Big(|\hat q_{1}(k)|^{p} + |\hat q_{2}(k)|^{p}  \Big)\, dk \le  \frac{C}{ \eps^{3d}}, 
\end{multline}
where $\xi_{1, \eps}$ and $\xi_{2, \eps}$ are given by \eqref{def-xi}.
It follows that 
\begin{equation*}
\|\nabla  (K_{\xi_{j, \eps}} * q_j)\|_{L^{2}(B_{r})} +   n \| K_{\xi_{j, \eps}} * q_j \|_{L^{2}(B_{r})} \le C_{r}/ \eps^{3d/p}, 
\end{equation*}
for all $r > 0$ and for $j=1, 2$.   By Proposition~\ref{pro1}, for $n$ large enough,  there exist $w_{j, \eps} \in H^1_{\loc}(\mR^d)$ ($j=1, 2$) such that
\begin{equation*}
w_{j, \eps} = K_{\xi_{j, \eps}} * (q_j + q_j w_{j, \eps}) \mbox{ in } \mR^d. 
\end{equation*}
and 
\begin{equation*}
\| \nabla w_{j, \eps}\|_{L^{2}(B_{r})} +s_{\eps} \| w_{j, \eps}\|_{L^{2}(B_{r})} \le C_{r}/ \eps^{3d}. 
\end{equation*}
The proof is complete. \proofend

\subsection{Proof of Theorem~\ref{thm1}}
Without loss of generality one may assume that $\Omega \subset B_{1/2}$.  Let   $g_{i, j}$,  $1 \le i, j \le 2$,  such that 
\begin{equation*}
g_{1, 1} + g_{1, 2} = g_{1} \mbox{ and } h_{1, 1} + h_{1, 2} = h_{1}. 
\end{equation*}
\begin{equation*}
g_{1, 1}, \; h_{1, 1} \mbox{ are smooth with compact support in $\Omega$}, 
\end{equation*}
\begin{equation*}
\|g_{1, 2}\|_{L^{\infty}} + \|h_{1, 2} \|_{L^{d}} \le 2 c,  
\end{equation*}
Extend  $g_{1,1}$ and $h_{1,1}$ smoothly in $\mR^d \setminus \Omega$ with compact support in $B_1$ and denote these extension by $G_{1,1}$ and $H_{1,1}$. 
Extend $g_{1, 2}, \, h_{1, 2}, \, g_2, \, h_2$ by 0 outside $\Omega$ and  denote these extensions by $G_{1, 2}, \, H_{1, 2}, \, G_2, \, H_2$. Define 
\begin{equation*}
G_1 = G_{1,1} + G_{1, 2} \quad \mbox{ and } \quad H_{1} = H_{1, 1} + H_{1, 2}. 
\end{equation*}
Extend $q_1$ and $q_2$ in $\mR^d$ by $ \dive G_1 + G_2 $ and $ \dive H_1 + H_2$ and still denote these extensions by  $q_1$ and $q_2$. 
Then $q_1$ and $q_2$ satisfy the assumptions of Proposition~\ref{pro2} since
$F( 1_{\Omega} ) \in L^r(\mR^d)$ 
{for $r > 2d/ (d+1)$ (see \cite[Theorem 1]{Lebelev}).} 
We claim that there exist $\sigma \in \mS^{d-1}$ and $u_{1, n}, \, u_{2, n} \in H^{1}_{loc}(\mR^{d})$ such that 
\begin{equation}\label{imp-0}
|\sigma - \sigma_{0}| \le \eps, 
\end{equation}
\begin{equation}\label{imp-1}
\Delta w_{j, n} + \xi_{j, n} \cdot \nabla w_{j, n} - q_{j} w_{j, n} = q_{j} \mbox{ in } \mR^{d}  \quad \mbox{ for } j=1, 2, 
\end{equation}
\begin{equation}\label{imp-2}
w_{j, n} \to 0 \mbox{ weakly in } H^{1}(B_{2}) \quad \mbox{ for } j=1, 2,
\end{equation}
for some $\xi_{1, n}, \xi_{2, n} \in \mC^{d}$ with $\xi_{j, n} \cdot \xi_{j, n} = 0$, $|\xi_{1, n} + \xi_{2, n} - \sigma | \to 0$ and $|\xi_{j, n} | \to \infty$ as $n \to \infty$.

\medskip
Indeed, by Proposition~\ref{pro2}, there exist $w_{j, n} \in H^{1}_{\loc}(\mR^{d})$ $(j=1, \, 2)$  such that 
\begin{equation*}
w_{j, n} = K_{\xi_{j, n}} * (q_{j}  + q_{j} w_{j, n}). 
\end{equation*}
Moreover, 
\begin{equation*}
\| \nabla w_{j, n}\|_{L^{2}(B_{r})} + n \| w_{j, n}\|_{L^{2}(B_{r})} \le C_{r}/ \eps^{3d/p},  
\end{equation*}
for some $C_{r} > 0$ which depends only on $d$, $g_{i}$, $h_{i}$ ($i=1, \, 2$), and $r$. Here 
\begin{equation*}
\xi_{1, n} =  s_{n} \sigma_{2, n}  - i s_{n} \sigma_{1, n}
\end{equation*}
and
\begin{equation*}
\xi_{2, n}=  s_{n} \Big( -  \frac{s_{\eps} \sigma_{2, n }}{ \sqrt{1 + s_{n}^{2}} }+ \frac{\sigma_{3, n}}{ \sqrt{1 + s_{n}^{2}}} \Big)  + i s_{n} \sigma_{1, n}. 
\end{equation*}
for some $s_{n} \in (n, 4n)$, and $\sigma_{1, n}, \sigma_{2, n}, \sigma_{3, n}$ such that
\begin{equation*}
\sigma_{1, n} \cdot \sigma_{2, n} = \sigma_{1,n} \cdot \sigma_{3, n} = \sigma_{2, n} \cdot \sigma_{3, n} = 0, 
\end{equation*}
\begin{equation*}
|\sigma_{3, n} - \sigma_{0} | \le \eps.
\end{equation*}
Without loss of generality one might assume that $\sigma_{3, n} \to \sigma$ for some $\sigma \in \mS^{d-1}$. Then 
\begin{equation*}
\xi_{1,n} + \xi_{2, n} =  \Big(s_n - \frac{s_n^2}{\sqrt{1 + s_n^2}} \Big) \sigma_{2,n} + \frac{s_n \sigma_{3, n}}{\sqrt{1 + s_n^2}}  \to \sigma,
\end{equation*}
and the claim is proved. 

\medskip
We now apply the complex geometric optics approach introduced by Sylvester and Uhlmann in \cite{SylvesterUhlmann87}. 
Define, for $j=1, 2$, 
\begin{equation*}
v_{j, n} = (1 + w_{j, n}) e^{\xi_{j, n} \cdot x/ 2}. 
\end{equation*}
Since $w_{j, n}$ satisfies \eqref{imp-1}, it follows that 
\begin{equation*}
\Delta v_{j, n} - q_{j} v_{j, n} = 0 \mbox{ in } \mR^{d}  \quad \mbox{ for } j=1, 2.
\end{equation*}
We derive from \eqref{imp-ident} that 
\begin{equation}\label{imp-3-0-1}
\int_{B_{2}} (q_{1} - q_{2})  (1 + w_{1, n}) (1 + w_{2, n}) e^{\sigma_{n} \cdot x / 2} = 0,  
\end{equation}
where 
\begin{equation}\label{imp-3-1}
\sigma_{n} = \xi_{1, n} + \xi_{2, n} \to \sigma \mbox{ as } n \to \infty. 
\end{equation}
A combination of  \eqref{imp-1}, \eqref{imp-3-0-1},  and \eqref{imp-3-1} yields
\begin{equation*}
\int_{B_{2}} (q_{1} - q_{2})  e^{\sigma \cdot x / 2} = 0. 
\end{equation*} 
Since $\sigma_{0} \in \mS^{d-1}$ and $\eps > 0$ are arbitrary, it follows that 
\begin{equation*}
\int_{B_{2}} (q_{1} - q_{2})  e^{\sigma_{0} \cdot x / 2} = 0 \mbox{ for all } \sigma_{0} \in \mS^{d-1}. 
\end{equation*} 
This implies 
\begin{equation*}
q_{1} = q_{2},
\end{equation*}
and the proof is complete. \proofend

\subsection{Proof of Corollary~\ref{cor1}}

Let $u_{i} \in H^{1}(\Omega)$ ($i=1, \, 2$) be a solution to the equation 
\begin{equation*}
\dive(\gamma_{i} \nabla u_{i}) = 0 \mbox{ in } \Omega. 
\end{equation*} 
Define 
\begin{equation*}
v_{i} = \gamma_{i}^{1/2} u_{i} \mbox{ in } \Omega. 
\end{equation*}
Then $v_{i} \in H^{1}(\Omega)$ is a solution to the equation 
\begin{equation*}
\Delta v_{i} - q_{i} v_{i} = 0 \mbox{ in } \Omega,  
\end{equation*}
where 
\begin{equation*}
q_{i} = \frac{\Delta \gamma_{i}^{1/2}}{\gamma_{i}^{1/2}} = \Delta t_{j} - |\nabla t_{j}|^{2} \mbox{ in } \Omega. 
\end{equation*} 
Here $t_{i}$ ($i=1, \, 2$) is given by
\begin{equation*}
t_{i} = \ln \gamma_{i}^{1/2} \mbox{ in } \Omega. 
\end{equation*}
Since $DtN_{\gamma_{1}} = DtN_{\gamma_{2}}$, it follows that 
\begin{equation*}
\int_{\Omega} (q_{1} - q_{2}) v_{1} v_{2} = 0, 
\end{equation*}
for all solutions $v_{i}$ ($i=1, \, 2$) to the equation \begin{equation*}
\Delta v_{i} - q_{i} v_{i} = 0 \mbox{ in } \Omega. 
\end{equation*}

Set 
\begin{equation*}
g_{i} = \nabla t_{i} \mbox{ and } h_{i} = - t_{i}^{2},
\end{equation*}
then $g_{i}$ and $h_{i}$ satisfy the assumptions of Theorem~\ref{thm1}. Applying Theorem~\ref{thm1}, we have
\begin{equation}\label{q1q2}
q_{1} = q_{2} \mbox{ in } \Omega. 
\end{equation}
This implies, 
\begin{equation*}
\Delta(t_1 - t_2) = |\nabla t_1|^2 - |\nabla t_2|^2 \in L^2(\Omega). 
\end{equation*}
Hence
\begin{equation}\label{UCP1}
\partial_{\eta} t_{1} = \partial_{\eta} t_{2} \mbox{ on } \partial \Omega, 
\end{equation}
We also have, by Proposition~\ref{Pro-boundary},  
\begin{equation}\label{UCP2}
t_{1} = t_{2} \mbox{ on } \partial \Omega, 
\end{equation}
We derive from \eqref{q1q2} and the definition of $q_{i}$ that
\begin{equation*}
\Delta (t_{1} - t_{2}) - \nabla t \cdot  \nabla (t_{1} - t_{2}) = 0 \mbox{ in } \Omega,
\end{equation*}
where $t = t_{1} + t_{2} \in W^{1, \infty}(\Omega)$. This implies $t_{1} = t_{2}$ by \eqref{UCP1}, \eqref{UCP2}, and the unique continuation principle. Therefore, the conclusion follows. \proofend

\section{Proof of Theorem~\ref{thm2} and Corollary~\ref{cor2}}
\label{lavinenachman}

\subsection{Construction of CGO solutions}

We begin this section with an estimate for the solution to the equation
\begin{equation*}
\Delta w + \xi \cdot \nabla w - q w = q \mbox{ in } \mR^{d}, 
\end{equation*}
for $q \in L^{d/2}(\mR^{d})$. 
This estimate will play an important role in the (new) proof of Theorem~\ref{thm2}.

\begin{proposition}\label{pro1-Nachmann} {Let  $d \ge 3$}, and let $\xi \in \mC^{d}$ with $\xi \cdot \xi = 0$, $q \in L^{d/2}(\mR^{d})$ with $\supp q \subset B_{1}$. For $|\xi|$ large enough,  there exists 
$w \in H^{1}_{\loc}(\mR^{d})$ such that
\begin{equation*}
w = K_{\xi}*( q + qw ). 
\end{equation*}
Moreover, 
\begin{equation}\label{est-u}
\lim_{|\xi| \to \infty}  \| w \|_{L^{\frac{2d}{d-2}}(B_{r})} =0, 
\end{equation}
\end{proposition}

\noindent{\bf Proof.} 
Let $f$ and $h$ be such that
\begin{equation*}
q = f + h,
\end{equation*}
where 
\begin{equation*}
f \mbox{ is smooth with support in $B_{1/2}$,  and } \| h\|_{L^{d/2}} \mbox{ is small}. 
\end{equation*}
Let $u_{0} = 0$  and consider the following iteration process: 
\begin{equation*}
 w_{n} = K_{\xi} * ( q + q w_{n-1} ) \mbox{ for } n \ge 1.
\end{equation*}
Define 
\begin{equation*}
w_{1, n}  =  K_{\xi} * (f + f w_{n-1})
\end{equation*}
and
\begin{equation*}
w_{2, n}  = K_{\xi} * (h + h w_{n-1}). 
\end{equation*}
Then
\begin{equation}\label{ident2}
w_{n} = w_{1, n} + w_{2,n}, 
\end{equation}
\begin{equation}\label{u1}
 w_{1, n+1} - w_{1, n}  =  K_{\xi} * \big[f(w_{n} - w_{n-1}) \big], 
\end{equation}
and
\begin{equation}\label{u2}
w_{2, n+1} - w_{2, n}  =   K_{\xi} * \big[ h (w_{n} - w_{n-1}) \big]. 
\end{equation}
Applying the generalized Sobolev's inequality  \cite[Theorem 2.1]{KRS} (see also \eqref{Sobolev}) and using \eqref{u2}, we have
\begin{align*}
\| w_{2, n+1} - w_{2, n}\|_{L^{\frac{2d}{d-2}}} \le & C \| \Delta (w_{2, n+1} - w_{2, n}) + \xi \cdot \nabla  (w_{2, n+1} - w_{2, n}) \|_{L^{\frac{2d}{d+2}}}  \nonumber \\[6pt]
\le & C \|h (w_{n} - w_{n-1}) \|_{L^{{\frac{2d}{d+2}}}};
\end{align*}
which yields
\begin{equation}\label{estimate-Sogge}
\| w_{2, n+1} - w_{2, n}\|_{L^{\frac{2d}{d-2}}} \le C \| h\|_{L^{d/2}} \|w_{n} - w_{n-1}\|_{L^{\frac{2d}{d-2}}(B_{1})}.
\end{equation}
We also have
\begin{equation*}
 \| w_{1, n+1} - w_{1, n} \|_{L^{\frac{2d}{d-2}}(B_{r})} \le  C_{r} \| w_{1, n+1} - w_{1, n} \|_{H^{1}(B_{r})}, 
\end{equation*}
and by \eqref{SU2}, 
\begin{align*}
\| w_{1, n+1} - w_{1, n} \|_{H^{1}(B_{r})} \le C_{r} \| f(w_{n} - w_{n-1}) \|_{L^{2}(B_{1})}. 
\end{align*}
This implies
\begin{equation}\label{interpolation2}
 \| w_{1, n+1} - w_{1, n} \|_{L^{\frac{2d}{d-2}}(B_{r})} \le  C_{r} \|f \|_{L^{\infty}} \| w_{n} - w_{n-1} \|_{L^{2}(B_{1})}. 
\end{equation}
A combination of \eqref{estimate-Sogge} and \eqref{interpolation2} yields
\begin{equation}\label{hahaha1}
\| w_{n+1} - w_{n}\|_{L^{\frac{2d}{d-2}}(B_{r})} 
\le C_{r} \Big( \| h\|_{L^{d/2}}\|w_{n} - w_{n-1}\|_{L^{\frac{2d}{d-2}}(B_{1})} \\[6pt]
+ \| f\|_{L^\infty} \|w_{n} - w_{n-1}\|_{L^{2}(B_{1})} \Big).
\end{equation}
On the other hand, by \eqref{SU1}, it follows from \eqref{u1} that 
\begin{align}\label{interpolation1}
 \| w_{1, n+1} - w_{1, n} \|_{L^{2}(B_{r})} \le & \frac{C_{r}}{|\xi|} \| f(w_{n} - w_{n-1}) \|_{L^{2}(B_{1})} \nonumber \\[6pt]
\le  & \frac{C_{r}}{|\xi|} \|f\|_{L^{\infty}} \|w_{n} - w_{n-1}\|_{L^{2}(B_{1})}.
\end{align}
A combination of \eqref{estimate-Sogge} and \eqref{interpolation1} implies
\begin{equation}\label{hahaha2}
\| w_{n+1} - w_{n}\|_{L^{2}(B_{r})} 
\le C_{r} \Big( \| h\|_{L^{d/2}}
+  \frac{\| f\|_{L^\infty}}{|\xi|} \Big) \|w_{n} - w_{n-1}\|_{L^{\frac{2d}{d-2}}(B_{1})}.
\end{equation}
From \eqref{hahaha1} and \eqref{hahaha2}, we obtain
\begin{multline}\label{tototo}
\| w_{n+1} - w_{n}\|_{L^{\frac{2d}{d-2}}(B_{r})} 
\le C_{r} \Big( \| h_1\|_{L^{d/2}}\|w_{n} - w_{n-1}\|_{L^{\frac{2d}{d-2}}(B_{1})} \\[6pt]
+ \| f_1\|_{L^{\infty}} \Big[\|h_2 \|_{L^{d/2}} + \frac{\| f_2\|_{L^\infty}}{|\xi|} \Big] \|w_{n-1} - w_{n-2}\|_{L^{\frac{2d}{d-2}}(B_{1})} \Big).
\end{multline}
Here $f_1, f_2$ and $h_1, h_2$ are such that
\begin{equation*}
q = f_1 + h_1 = f_2 + h_2,
\end{equation*}
where 
\begin{equation*}
f_1, f_2 \mbox{ are smooth with support in $B_{1}$,  and } \| h_1\|_{L^{d/2}}, \| h_2\|_{L^{d/2}}  \mbox{ are small}. 
\end{equation*}
Appropriate choice of $f_1, f_2$ and $h_1, h_2$ implies that $w_n$ converges to $w$ in $H^1_{\loc}(\mR^d)$ and 
\begin{equation*}
w = K_{\xi}*( q + qw ). 
\end{equation*}

\medskip
We next prove \eqref{est-u}. 
By the same arguments used to obtain~\eqref{hahaha1} and \eqref{hahaha2}, we have
\begin{equation}\label{ttt1}
\| w - K_\xi * q \|_{L^{\frac{2d}{d-2}}(B_{r})} 
\le C_{r} \Big( \| h_1\|_{L^{d/2}}\|w\|_{L^{\frac{2d}{d-2}}(B_{1})} 
+ \| f_1\|_{L^{\infty}} \|w\|_{L^{2}(B_1)} \Big)
\end{equation}
and 
\begin{equation}\label{ttt2}
\| w - K_\xi * q \|_{L^{2}(B_{r})} 
\le C_{r} \Big( \| h_2\|_{L^{d/2}}\|w\|_{L^{\frac{2d}{d-2}}(B_{1})} 
+ \frac{\| f_2\|_{L^{\infty}}}{|\xi|} \|w\|_{L^{2}(B_1)} \Big).
\end{equation}
We claim that 
\begin{equation}\label{ttt3}
\lim_{|\xi| \to 0}\| K_\xi * q\|_{L^{\frac{2d}{d-2}}(B_1)} = 0. 
\end{equation}
Admitting \eqref{ttt3}, we will prove \eqref{est-u}. In fact, from \eqref{ttt2}, and \eqref{ttt3}, we have $\| w \|_{L^2(B_r)} \to 0$ as $|\xi| \to \infty$. This implies, by \eqref{ttt3}, $\| w \|_{L^\frac{2d}{d-2}(B_r)} \to 0$ as $|\xi| \to \infty$. 

\medskip
It remains to  prove \eqref{ttt3}. 
Let $q_1, q_2 \in L^{d/2}$ with the support  in $B_{1}$ be such that 
\begin{equation*}
q = q_1 + q_2, \quad  \| q_1 \|_{L^{d/2}} \mbox{ is small},  \quad \mbox{ and } \quad  q_2 \mbox{ is smooth}. 
\end{equation*}
We have 
\begin{equation*}
\| K_\xi * q\|_{L^{\frac{2d}{d-2}}(B_1)} \le \| K_\xi * q_1\|_{L^{\frac{2d}{d-2}}(B_1)} + \| K_\xi * q_2\|_{L^{\frac{2d}{d-2}}(B_1)}, 
\end{equation*}
and,  by the generalized Sobolev inequality,  
\begin{equation*}
\| K_\xi * q_1\|_{L^{\frac{2d}{d-2}}}   \le C \|q_1 \|_{L^{\frac{d}{2}}}, 
\end{equation*}
and by \eqref{SU1}, 
\begin{equation*}
\| K_\xi * q_2\|_{L^{\frac{2d}{d-2}}(B_1)} \le  \frac{C}{|\xi|} \| q_2\|_{C^2}. 
\end{equation*}
By an appropriate choice of $q_1$ and $q_2$, it follows that 
\begin{equation*}
\lim_{|\xi| \to 0}\| K_\xi * q\|_{L^{\frac{2d}{d-2}}(B_1)} = 0; 
\end{equation*}
claim \eqref{ttt3} is proved.  The proof is complete. \proofend

\subsection{Proof of Theorem~\ref{thm2}}

The proof is standard after Proposition~\ref{pro1-Nachmann}. For the convenience of the reader, we present the proof.  Without loss of generality one may assume that $\Omega \subset B_{1/2}$.  Extend $q_1$ and $q_2$ by $0$ in $\mR^d \setminus \Omega$ and still denote these extensions by $q_1$ and $q_2$.  Let  $\sigma_{1}, \sigma_{2}, \sigma_{3} \in \mS^{d-1}$ be such that 
\begin{equation*}
\sigma_{1} \cdot \sigma_{2} = \sigma_{1} \cdot \sigma_{3} = \sigma_{2} \cdot \sigma_{3} = 0. 
\end{equation*}
Set 
\begin{equation*}
\xi_{1, n} =  n \sigma_{2}  - i n \sigma_{1} \mbox{ and }
\xi_{2, n}=  n \Big( - \frac{n\sigma_{2}}{ \sqrt{1 + n^{2}}} + \frac{ \sigma_{3}}{ \sqrt{1 + n^{2}}} \Big)  + i n \sigma_{1}. 
\end{equation*}
For $n$ large enough, by Proposition~\ref{pro1-Nachmann}, there exist $w_{j, n} \in H^1_{\loc}$ ($j=1, 2$) such that 
\begin{equation*}
w_{j, n} = K_{\xi_{j, n}}*( q_j + q_j u_{j, n}). 
\end{equation*}
Moreover, 
\begin{equation}\label{est-u-i}
\lim_{n \to \infty}  \| w_{j, n} \|_{L^{\frac{2d}{d-2}}(B_{r})} =0 \quad \mbox{ for } j =1,2. 
\end{equation}
Define, for $j=1, 2$, 
\begin{equation*}
v_{j, n} = (1 + w_{j, n}) e^{\xi_{j, n} \cdot x/ 2}, 
\end{equation*}
Then 
\begin{equation*}
\Delta v_{j, n} + q_{j} v_{j, n} = 0 \mbox{ in } \mR^{d}  \quad \mbox{ for } j=1, 2.
\end{equation*}
We derive from \eqref{imp-ident} that 
\begin{equation}\label{imp-3-0}
\int_{B_{2}} (q_{1} - q_{2})  (1 + w_{1, n}) (1 + w_{2, n}) e^{\sigma_{s} \cdot x / 2} = 0,  
\end{equation}
where 
\begin{equation}\label{imp-3}
\sigma_{s} = \xi_{1, n} + \xi_{2, n} \to \sigma_{3} \mbox{ as } s \to \infty. 
\end{equation}
A combination of  \eqref{est-u-i}, \eqref{imp-3-0},  and \eqref{imp-3} yields
\begin{equation*}
\int_{B_{2}} (q_{1} - q_{2})  e^{\sigma \cdot x / 2} = 0. 
\end{equation*} 
Since $\sigma_{3}  \in \mS^{d-1}$ is arbitrary, it follows that 
\begin{equation*}
q_{1} = q_{2}, 
\end{equation*}
and the proof is complete. \proofend

\subsection{Proof of Corollary~\ref{cor2}}

The proof is similar to the one of Corollary~\ref{cor1}. The details are left to the reader. \proofend 
 
\section{Uniqueness of Calderon's problem \\ for  conductivities of class $W^{s, 3/s}$ for $s > 3/2$ in $3d$}
\label{newclass}

\subsection{Construction of CGO solutions}

We begin this section with 
\begin{lemma}\label{lem-new} Let  $\xi \in \mC^{3}$ with $|\xi| > 2$ and $\xi \cdot \xi =0$, $v \in H^{1}_{\loc}(\mR^{3})$ and $q \in H^{-1/2}(\mR^3)$ with $\supp q  \subset B_{1}$. 
Define 
\begin{equation*}
W = K_{\xi} *(qV). 
\end{equation*}
We have
\begin{equation*}
\| W \|_{H^{1}(B_{r})}^{2} \le C_{r} \|V\|_{H^{1}(B_2)}^{2} \cdot   E(q, \xi), 
\end{equation*}
where
\begin{multline*}
E(q, \xi) = \int_{\mR^{3}} |\hat q(\eta)|^{2} \int_{4 |\xi| \ge \dist(k, \Gamma_{\xi}) \ge |k|/ |\xi|} \frac{|k|^{2} |\hat K_{\xi}(k)|^{2}}{|k - \eta|^{2}}  \, dk \, d \eta \\[6pt]
+  \int_{\mR^{3}} |\tilde q(\gamma)|^{2} \int_{\dist(\eta, \Gamma_{\xi}) \le |\eta|/ |\xi|} \frac{1}{|\eta - \gamma|^{2}} \, d \eta \, d \gamma  + 
\| q\|_{H^{-1/2}}^2, 
\end{multline*}
where
\begin{equation*}
\widetilde q(k) : = \sup_{\eta \in B_{4}(k)} |\hat q(\eta)|. 
\end{equation*}
\end{lemma}

\noindent{\bf Proof.} Without loss of generality, one may assume that $\supp V \subset B_{3/2}$ and $r > 1$. Set 
\begin{equation}\label{def-f}
f = q V,
\end{equation}
then 
\begin{equation*}
W = K_{\xi} * f. 
\end{equation*}

We first prove 
\begin{equation}\label{partL2}
\| W \|_{L^2(B_r)} \le C_r\|V \|_{H^1} \| q\|_{H^{-1/2}}. 
\end{equation}
Applying \eqref{SU2}, we have 
\begin{equation}\label{uL2}
\| W\|_{L^2(B_r)} \le C_r \|f\|_{H^{-1}}. 
\end{equation}
On the other hand, 
\begin{equation}\label{fH-1-1}
\|f\|_{H^{-1}}^2 \le \int_{\mR^3} \frac{1}{|k|^2 + 1} \Big| \int_{\mR^3} |\hat q(k - \eta)| |\hat V (\eta)| \, d \eta \Big|^2 \, d k. 
\end{equation}
Since 
\begin{equation*}
\int_{\mR^3} |\hat q(k - \eta)| |\hat V (\eta)| \, d \eta = \int_{|\eta| \le |k|/2 } |\hat q(k - \eta)| |\hat V (\eta)| \, d \eta + \int_{|\eta | >  |k|/2} |\hat q(k - \eta)| |\hat V (\eta)| \, d \eta, 
\end{equation*}
it follows that, by H\"older's inequality,  
\begin{multline}\label{fH-1-2}
\frac{1}{2}\int_{\mR^3} \frac{1}{|k|^2 + 1} \Big| \int_{\mR^3} |\hat q(k - \eta)| |\hat V (\eta)| \, d \eta \Big|^2 \, d k \\[6pt] 
\le \int_{\mR^3}  \int_{|\eta| \le |k|/2} \frac{|\hat q(k - \eta)|^2}{(k^2 + 1)(|\eta|^2 + 1)}  \, d \eta \int_{|\eta| \le  |k|/2}  |\hat V (\eta)|^2 (|\eta|^2 +1 )\, d \eta \, d k \\[6pt]
 + \int_{\mR^3}  \int_{|\eta| \ge  |k|/2} \frac{|\hat q(k - \eta)|^2}{(|k - \eta|^2 + 1)^{1/2}} \, d \eta  \int_{|\eta| \ge  |k|/2}  \frac{ |\hat V (\eta)|^2 (|k -\eta|^2 +1 )^{1/2}}{|k|^2 + 1} \, d \eta \, d k. 
\end{multline}
We have, since $|k - \eta| \le |k|/2$ implies $2 |\eta| \ge |k| \ge 2 |\eta|/3$,  
\begin{multline}\label{fH-1-3}
 \int_{\mR^3}  \int_{|\eta| \le |k|/2} \frac{|\hat q(k - \eta)|^2}{(k^2 + 1)(|\eta|^2 + 1)}  \, d \eta \, d k 
 =\int_{\mR^3}  \int_{|k - \eta| \le |k|/2} \frac{|\hat q(\eta)|^2}{(k^2 + 1)(|k - \eta|^2 + 1)}   \, d k \, d \eta \\[6pt]
 \le \int_{\mR^3}\frac{|\hat q(\eta)|^2}{(1 + |\eta|^2)^{1/2}} \int_{2 |\eta| \ge |k| \ge 2 |\eta|/3} \frac{1}{(|k - \eta|^2 +1)(1 + |\eta|^2)^{1/2}} \le C
  \| q\|_{H^{-1/2}}^2
\end{multline}
and 
\begin{equation}\label{fH-1-4}
\int_{\mR^3}  \int_{|\eta| \ge  |k|/2}  \frac{ |\hat V (\eta)|^2 (|k -\eta|^2 +1 )^{1/2}}{|k|^2 + 1} \, d \eta \, d k \le C \| V\|_{H^1}^2. 
\end{equation}
Using \eqref{fH-1-1}, \eqref{fH-1-2}, \eqref{fH-1-3}, and \eqref{fH-1-4}, we derive from \eqref{fH-1-1} that 
\begin{equation}\label{fH-1}
\|f\|_{H^{-1}} \le C \|V \|_{H^1} \| q\|_{H^{-1/2}}. 
\end{equation}
A combination of  \eqref{uL2} and \eqref{fH-1} yields \eqref{partL2}. 

\medskip
It remains to establish the key estimate
\begin{equation}\label{conclusion-1}
\|\nabla W\|_{L^2(B_r)} \le C_r \|V\|_{H^{1}} \cdot E(q, \xi). 
\end{equation}
Set
\begin{equation*}
\Gamma_{\xi} := \{k \in \mR^{3} ; \; - |k|^{2}  + i \xi \cdot k =0 \}. 
\end{equation*}
Define $K_{1, \xi}$, $K_{2, \xi}$, and $K_{3, \xi}$ as follows
\begin{equation*}
\hat K_{1, \xi} (k) = 
\left\{\begin{array}{cl}
\hat K_{\xi}(k) &  \mbox{ if } 4 |\xi| \ge \dist(k, \Gamma_\xi) > |k|/|\xi| ,  \\[6pt]
0 &  \mbox{ otherwise}, 
\end{array} \right. 
\end{equation*}
\begin{equation*}
\hat K_{2, \xi} (k) = \left\{\begin{array}{cl}
\hat K_{\xi}(k) &  \mbox{ if }  \dist(k, \Gamma_\xi) \le  |k|/|\xi| ,  \\[6pt]
0 &  \mbox{ otherwise}, 
\end{array} \right. 
\end{equation*}
and
\begin{equation*}
\hat K_{3, \xi} (k) = \left\{\begin{array}{cl}
\hat K_{\xi}(k) &  \mbox{ if }  \dist(k, \Gamma_\xi) > 4 |\xi| ,  \\[6pt]
0 &  \mbox{ otherwise}.
\end{array} \right. 
\end{equation*}
Then
\begin{equation}\label{b1-1} \| \nabla (K_{\xi} * f) \|_{L^{2}(B_{r})}  
 \le \| \nabla (K_{1, \xi} * f) \|_{L^{2}(B_{r})}  + \| \nabla (K_{2, \xi} * f) \|_{L^{2}(B_{r})} +  \| \nabla (K_{3, \xi} * f) \|_{L^{2}(B_{r})}  . 
\end{equation}
Since 
\begin{equation*}
|\hat K_{\xi}(k)| \le \frac{1}{ |k|^2} \mbox{ for } \dist(k, \Gamma_\xi)  \ge 4 |\xi|, 
\end{equation*}
it follows that 
\begin{equation}\label{K31}
 \| \nabla (K_{3, \xi} * f) \|_{L^{2}(B_{r})}  \le  \| \nabla (K_{3, \xi} * f) \|_{L^{2}(\mR^{3})} \le C \| f\|_{H^{-1}}. 
 \end{equation}
A combination of \eqref{fH-1} and \eqref{K31}  yields 
\begin{equation}\label{K3}
 \| \nabla (K_{3, \xi} * f) \|_{L^{2}(B_{r})}  \le C \|V \|_{H^1} \| q\|_{H^{-1/2}}. 
 \end{equation}
We next estimate the first two terms in the RHS of \eqref{b1-1}. We start with $\|\nabla (K_{1, \xi} * f) \|_{L^2(B_r)}$.  
Since 
\begin{equation*}
 \| \nabla (K_{1, \xi} * f) \|_{L^{2}(B_{r})}^{2} \le \| \nabla (K_{1, \xi} * f) \|_{L^{2}(\mR^{3})}^{2} ,
\end{equation*}
it follows from Plancherel's theorem that 
\begin{equation}\label{b22-1}
 \| \nabla (K_{1, \xi} * f) \|_{L^{2}(B_{r})}^{2} \le   C \int_{4|\xi| \ge |\hat K_{\xi} (k)| \ge k/|\xi|} |\hat f(k)|^{2} |k|^{2} |\hat K_{\xi} (k)|^{2} \, dk. 
\end{equation}
From \eqref{def-f}, we have 
\begin{equation*}
\hat f(k) = \int_{\mR^{3}} \hat q(\eta) \hat V (k - \eta) \, d \eta. 
\end{equation*}
Applying H\"older's inequality, we obtain 
\begin{equation}\label{b22-2}
|\hat f(k)|^{2} \le  \int_{\mR^{3}} \frac{|\hat q(\eta)|^{2}}{|k - \eta|^{2}} \, d \eta  \int_{\mR^{3}} |(k - \eta) \hat V (k - \eta)|^{2} \, d \eta. 
\end{equation}
A combination of \eqref{b22-1} and \eqref{b22-2} yields 
\begin{multline}\label{K1}
\int_{4|\xi| \ge |\hat K_{\xi} (k)| \ge k/|\xi|} |\hat f(k)|^{2} |k|^{2} |\hat K_{\xi} (k)|^{2} \, dk\\[6pt]
 \le  C \| \nabla V \|_{L^{2}}^{2}   \int_{\mR^{3}} |\hat q(\eta)|^{2} \int_{4 |\xi| \ge |\hat K_{\xi} (k)| \ge |k|/ |\xi|} \frac{|k|^{2} |\hat K_\xi(k)|^{2}}{|k - \eta|^{2} } \, dk \, d \eta. 
\end{multline}
We next estimate $\| \nabla (K_{2, \xi} *f) \|_{L^2(B_r)}$. Fix 
\begin{equation*}
\varphi \in C^{\infty}_{0}(\mR^{3}) \mbox{ with } \varphi = 1 \mbox{ in } B_{2 r}.
\end{equation*} 
Define 
\begin{equation}\label{def-f1-1}
\widetilde f(k) = \sup_{\eta \in B_{4}(k)} |\hat f(\eta)|, 
\end{equation}
and
\begin{equation}\label{def-phi1-1}
\widetilde \varphi (k) = \sup_{\eta \in B_{4}(k)} |\hat \varphi (\eta) | .
\end{equation}
Since 
\begin{equation*}
 |\hat f|* |\hat \varphi| (\eta) = \int_{\mR^{d}} |\hat f (\zeta)| |\hat \varphi(\eta -\zeta)| \, d \zeta, 
\end{equation*}
and $f = f \varphi$,  
it follows from the definition of $\widetilde f$ \eqref{def-f1-1} and  $\widetilde \varphi$ \eqref{def-phi1-1} that 
\begin{equation*}
\widetilde f \le |\hat f|* \widetilde \varphi. 
\end{equation*}
Since 
\begin{equation*}
\|\nabla( K_{2, \xi} * f) \|_{L^{2}(B_{r})}^{2} \le \|\nabla (\varphi  \cdot [K_{2, \xi} * f]) \|_{L^{2}(\mR^{3})}^{2}, 
\end{equation*}
it follows that 
\begin{equation}\label{to0-1}
\|\nabla( K_{2, \xi} * f) \|_{L^{2}(B_{r})}^{2}  \le  \int_{\mR^{d}} |k|^2 \left| \int_{\dist(\eta, \Gamma_{\xi}) \le |\eta|/ |\xi|} |\hat \varphi(k- \eta)| \cdot  |\hat K_{\xi} (\eta)| \cdot |\hat f(\eta)| \, d \eta \right|^{2} \, dk. 
\end{equation}
Using the fact that $\hat K_\xi(\eta) \le  C / \big( |\xi| \dist (\eta, \Gamma_\xi)\big)$ for $|\eta| \le 2 |\xi|$ and 
\begin{equation*}
\int_{|x| \le 1} \frac{1}{|x_{1}| + |x_{2}|} \, dx < + \infty, 
\end{equation*}
as in \eqref{to1}, we obtain
\begin{equation}\label{to1-1}
\int_{\dist(\eta, \Gamma_{\xi}) \le 1} |\hat \varphi(k- \eta)| \cdot  |\hat K_{\xi} (\eta)| \cdot |\hat f(\eta)| \, d \eta \le  \frac{C}{|\xi|}\int_{\dist(\eta, \Gamma_{\xi}) \le |\eta|/ |\xi |}  \widetilde \varphi(k- \eta) \cdot \widetilde f(\eta) \, d \eta. 
\end{equation}
Applying H\"older's inequality, we derive from \eqref{to0-1} and \eqref{to1-1} that  
\begin{equation}\label{b3-1-1}
\|\nabla( K_{2, \xi} * f ) \|_{L^{2}(B_{r})}^{2} \le  \frac{C}{|\xi|^{2}} \int_{\mR^{3}} \int_{\dist(\eta, \Gamma_{\xi}) \le |\eta|/ |\xi|}  |k|^{2} \tilde \varphi(k - \eta) |\widetilde f (\eta)|^{2} \,  d \eta \, dk. 
\end{equation}
Since $|k|^{2} \le C (|k - \eta|^{2} + |\eta|^{2})$ and $\widetilde \varphi$ decays fast at infinity, it follows from \eqref{def-phi1-1} and \eqref{b3-1-1} that 
\begin{equation}\label{step1-K2}
\|\nabla( K_{2, \xi} * f ) \|_{L^{2}(B_{r})}^{2} \le  C \int_{\dist(\eta, \Gamma_{\xi}) \le |\eta|/ |\xi|}  |\widetilde f (\eta)|^{2} \,  d \eta. 
\end{equation}
Since 
\begin{equation*}
\tilde f (\eta) \le \tilde q * |\hat V| (\eta), 
\end{equation*}
it follows that 
\begin{equation*}
\int_{\dist(\eta, \gamma_\xi) \le |\eta|/ |\xi|} |\tilde f (\eta)|^2 \, d \eta  \le \int_{\dist(\eta, \gamma_\xi) \le |\eta|/ s} \left| \int_{\mR^3} \frac{|\tilde q(\gamma)|}{|\eta - \gamma|}  |\eta - \gamma | |\hat V (\eta - \gamma)| \, d \gamma \right|^2 \, d \eta. 
\end{equation*}
Using H\"older's inequality, we obtain
\begin{equation}\label{step2-K2}
\int_{\dist(\eta, \Gamma_{\xi}) \le |\eta|/ |\xi|}  |\widetilde f (\eta)|^{2} \,  d \eta \le C \| \nabla  V\|_{L^{2}(\mR^{3})}^{2} \int_{\mR^{3}} |\tilde q(\gamma)|^{2} \int_{\dist(\eta, \Gamma_{\xi}) \le |\eta|/ |\xi|} \frac{1}{|\eta - \gamma|^{2}} \, d \eta \,  d \gamma. 
\end{equation}
A combination of \eqref{step1-K2} and \eqref{step2-K2} yields 
\begin{equation}\label{K2}
\|\nabla( K_{2, \xi} * f ) \|_{L^{2}(B_{r})}^{2} \le  C \| \nabla  V\|_{L^{2}(\mR^{3})}^{2} \int_{\mR^{3}} |\tilde q(\gamma)|^{2} \int_{\dist(\eta, \Gamma_{\xi}) \le |\eta|/ |\xi|} \frac{1}{|\eta - \gamma|^{2}} \, d \eta \,  d \gamma.
\end{equation}

We derive from \eqref{K3}, \eqref{K1}, and \eqref{K2} that \eqref{conclusion-1} holds. The proof is complete. \proofend

\medskip 
To use Lemma~\ref{lem-new}, we need to choose $\xi$ such that $E(q, \xi)$ remains bounded. This can be done using the following 
 average estimate for $E(q, \xi)$ whose proof is in the spirit of the one of Lemma~\ref{lem-average1} and is presented in the appendix. 

\begin{lemma}\label{lem-average-11} 
 Let $d = 3$ and  $R > 10$. We have 
\begin{equation}\label{avrg}
{1\over R} \int_{R/2}^{2R}
\int_{\mS^{d-1}} \int_{\mS^{d-1}_{\sigma_{1}}} E(q, s \sigma_{2} - i s \sigma_{1})\, d \sigma_{2} \, d \sigma_{1} \, ds  \le C \int_{\mR^{3}} |\hat q(\eta)|^{2}   \min \Big\{\frac{\ln R}{R}, \frac{R\ln R}{|\eta|^{2}} \Big\}\, d \eta
\end{equation} 
and
\begin{multline}\label{avrg2}
{1\over R} \int_{R/2}^{2R} \int_{\sigma_{1} \in \mS^{2}} \int_{\sigma_{2} \in \mS^{2}_{\sigma_{1}}} \int_{\sigma_{3} \in \mS^{2}_{\sigma_{1}, \sigma_{2}}} E \Big( q, \frac{s^2\sigma_{2}}{ \sqrt{1 + s^{2}}} + \frac{s \sigma_{3}}{ \sqrt{1 + s^{2}}}  - i s \sigma_{1} \Big) \, d \sigma_{3}\, d \sigma_{2} \, d \sigma_{1} \, ds\\
 \le C \int_{\mR^{3}} |\hat q(\eta)|^{2}   \min \Big\{\frac{\ln R}{R}, \frac{R \ln R}{|\eta|^{2}} \Big\} \, d \eta. 
\end{multline}
\end{lemma}

We recall that, by \eqref{S1} and \eqref{S12},  
\begin{equation*}
\mS^{2}_{\sigma_{1}} := \big\{ \sigma \in \mS^{2}; \; \sigma \cdot \sigma_{1} = 0 \big\}
 \end{equation*}
 and 
\begin{equation*}
\mS^{2}_{\sigma_{1}, \sigma_{2}} := \big\{ \sigma \in \mS^{2}; \; \sigma \cdot \sigma_{1} = 0 \mbox{ and } \sigma \cdot \sigma_{2} = 0 \big\}. 
 \end{equation*} 

\medskip
We will show that the RHS of \eqref{avrg} will behave like $\| q\|_{H^{-1/2}}$ for appropriate choice of $s$. For this end, we need the following lemma.

\begin{lemma}\label{lem-sequence} Let $(a_{n})$ be a non-negative sequence. Define 
\begin{equation*}
b_{n} = \sum_{l=1}^{n} 2^{l - n} a_{l}. 
\end{equation*}
Assume that $S = \sum_{1}^{\infty} a_{n} < + \infty$, then 
\begin{equation*}
\liminf_{n \to \infty} n b_{n} = 0. 
\end{equation*}
\end{lemma}

\noindent{\bf Proof. } The conclusion is a consequence of the following facts: 
\begin{equation*}
\sum_{n=1}^{\infty} b_{n} \le c \sum_{1}^{\infty} a_{n} < + \infty
\end{equation*}
for some positive constant $c$, and 
\begin{equation*}
\liminf_{n \to \infty} n b_{n} = 0, 
\end{equation*}
if 
\begin{equation*}
\sum_{n=1}^{\infty} b_{n}  < + \infty. 
\end{equation*}
\proofend

Applying Lemmas~\ref{lem-new}, \ref{lem-average-11}, and \ref{lem-sequence}, we can obtain the following result which is a variant of Propositions~\ref{pro1} and \ref{pro2} in this setting. 

\begin{proposition} \label{pro1-new} Let $q_{1}, q_{2} \in H^{-1/2}(\mR^{3})$ with support in $B_{1}$, and  $\sigma_{1}, \sigma_{2}, \sigma_{3} \in \mS^{2}$ be such that 
\begin{equation*}
\sigma_{1} \cdot \sigma_{2} = \sigma_{1} \cdot \sigma_{3} = \sigma_{2} \cdot \sigma_{3} =0. 
\end{equation*}
For any $\eps > 0$, there exist a sequence $s_{n} \to \infty$, $\sigma_{1, n}, \sigma_{2, n}, \sigma_{3, n} \in \mS^{2}$ and  $u_{i, n} \in H^{1}_{\loc}(\mR^{3})$ such that 
\begin{equation*}
\sigma_{1, n} \cdot \sigma_{2, n} = \sigma_{1, n} \cdot \sigma_{3, n} = \sigma_{2, n} \cdot \sigma_{3, n} =0, 
\end{equation*}
\begin{equation*}
|\sigma_{j, n} - \sigma_{j}| \le \eps \mbox{ for } j =1, 2, 3, 
\end{equation*}
and 
\begin{equation*}
w_{j, n} = K_{\xi_{j, n}} * (q_{j} + q_{j} w_{j, n}) \mbox{ for } j =1, 2. 
\end{equation*}
Here
\begin{equation*}
\xi_{1, n}  =  s_{n} \sigma_{2, n}  - i s \sigma_{1, n} \quad \mbox{ and } \quad
\xi_{2, n} =   - \frac{s_{n}^2 \sigma_{2, n}}{ \sqrt{1 + s_{n}^{2}}} + \frac{s_n \sigma_{3, n}}{ \sqrt{1 + s_{n}^{2}} }  + i s_{n} \sigma_{1, n}. 
\end{equation*}
Moreover, 
\begin{equation*}
\| w_{j, n}\|_{H^{1}(B_{r})} \le C \eps  \| K_{\xi_{j, n}} * q_{j}  \|_{H^1(B_1)} \quad \mbox{ for } j =1, 2,  
\end{equation*}
and for large $n$. 
\end{proposition}

\noindent{\bf Proof.} For $\eps > 0$, let $q_{j,1} \in C^\infty(\mR^3)$ and $q_{j, 2} \in C^\infty(\mR^3)$ with supports in $B_1$ be such that, for $j=1, 2$,  
\begin{equation*}
q_{j, 1} + q_{j, 2} = q_j, 
\end{equation*}
and 
\begin{equation}\label{q2}
\| q_{j, 2}\|_{H^{-1/2}} \le \eps^4. 
\end{equation} 
Define 
\begin{equation*}
a_{j, n} = \int_{2^{n} \le |k | \le 2^{n+1}} \frac{|\hat q_{j, 2}(k)|^{2}}{|k|} \, d k.  
\end{equation*}
It is clear that 
\begin{equation*}
\sum_{n=1}^{\infty} a_{j, n} \le \| q_{j, 2}\|_{H^{-1/2}}^{2}. 
\end{equation*}
Set 
\begin{equation*}
b_{j, n} := \sum_{l=1}^{n} 2^{l - n} a_{j, l} \sim \int_{2^{1} \le |k| \le 2^{n+1}} \frac{|\hat q_{j, 2}(k)|^{2}}{2^{n+1}} \, dk. 
\end{equation*}
By Lemma~\ref{lem-sequence}, there exists $n_{k} \to \infty$ such that 
\begin{equation*}
n_{k} b_{1, n_{k}} + n_{k} b_{2, n_{k}} \le c  \big( \| q_{1, 2}\|_{H^{-1/2}}^{2} + \| q_{2, 2}\|_{H^{-1/2}}^{2} ).  
\end{equation*}
Applying Lemma~\ref{lem-average-11}, there exist $\sigma_{1, k}, \sigma_{2, k}, \sigma_{3, k} \in \mS^{2}$ such that 
\begin{equation*}
\sigma_{1, k} \cdot \sigma_{2, k} = \sigma_{1, k} \cdot \sigma_{3, k} = \sigma_{2, k} \cdot \sigma_{3, k} =0, 
\end{equation*}
\begin{equation*}
|\sigma_{j, k} - \sigma_{j}| \le \eps \mbox{ for } j =1, 2, 3, 
\end{equation*}
and, by \eqref{q2},  
\begin{equation}\label{avrg-1-2}
E(q_{1, 2},  \xi_{1, k}) + E(q_{2, 2},  \xi_{2, k}) \le C \eps \mbox{ for large } k. 
\end{equation}
Here 
\begin{equation*}
\xi_{1, k}  = s_{k} \sigma_{2, k} - i s_{k} \sigma_{1, k} \quad \mbox{ and } \quad 
\xi_{2, k}  = -  \frac{ s_{k}^2 \sigma_{2, k} }{ \sqrt{1 + s_{k}^{2}} } + \frac{s_k \sigma_{3, k}} {\sqrt{1 + s_{k}^{2}} }   + i s_{k} \sigma_{1, k}. 
\end{equation*}
Let $w_{j, k, 0} = 0$ ($j=1, 2$) and consider the following iteration process: 
\begin{equation*}
 w_{j, k, n} = K_{\xi_{j, k}} * ( q_{j} + q_{j} w_{j, k, n-1} ) \mbox{ for } n \ge 1.  
\end{equation*}
Then, for $n \ge 1$ and $j=1, 2$,  
\begin{align*}
 w_{j, k,  n+1} - w_{j, k, n} = &K_{\xi_{j, k}} * (q_{j} [w_{j, k, n} - w_{j, k, n-1}] )\\[6pt]
  = & K_{\xi_{j, k}} * (q_{j, 1}[w_{j, k, n} - w_{j, k, n-1}] ) + K_{\xi_{j, k}} * (q_{j, 2} [w_{j, k, n} - w_{j, k, n-1}]). 
\end{align*}
Applying \eqref{SU1} for the first part and  applying Lemma~\ref{lem-new} for the second part, we have 
\begin{equation*}
\|w_{j, k, n+1}  - w_{j, k, n}\|_{H^{1}(B_{r})}^{2} \le C_{r} \Big( E(q_{j, 2}, \xi_{j, k})^{1/2} + \frac{\|q_{j, 1} \|_{C^2}}{|\xi_{j, k}|} \Big)\| w_{j, k, n}  - w_{j, k, n-1} \|_{H^{1}(B_{1})}^{2}. 
\end{equation*}
This implies 
\begin{equation}\label{haha1}
\sum_{m=1}^n \| w_{j, k, m+1} - w_{j, k, m}\|_{H^1(B_{r})} \le c(r, f, h) \sum_{m=0}^{n-1} \| w_{j, k, m} - w_{j, k,  m - 1}\|_{H^1(B_{1})}, 
\end{equation}
where $c(r, f, h) = C_r \Big( E(q_{i, 2}, \xi_{i, k})^{1/2} + \frac{\| q_{j, 1}\|_{C^2}}{|\xi_{i, k}|} \Big)$. 
From \eqref{avrg-1-2}, for large $k$, we derive that $c(2, f, h) < 1/2$. For such a large $k$, 
we have 
\begin{equation*}
\sum_{m=1}^n \| w_{j, k,  m+1} - w_{j, k, m}\|_{H^1(B_{2})} \le 2 c(r, f, h) \| w_{j, k,  1} - w_{j, k, 0}\|_{H^1(B_{1})} = 2 c(2, f, h) \| w_{j, k, 1} \|_{H^1(B_{1})}.
\end{equation*}
It follows that  there exist $w_{j} \in H^{1}_{\loc}(\mR^{3})$ such that 
\begin{equation*}
w_{j, k} = K_{\xi_{j, k}} * (q_{j, k} + q_{j} w_{j, k}), 
\end{equation*}
and, by \eqref{avrg-1-2} and \eqref{haha1},  
\begin{equation*}
\|w_{j, k} \|_{H^1(B_r)} \le C \eps \| w_{j, k, 1} \|_{H^1(B_{1})}, 
\end{equation*}
for large $k$ large enough.  The proof is complete. \proofend

\subsection{Proof of Theorem~\ref{thm3}}

Theorem~\ref{thm3} is a consequence of Proposition~\ref{pro1-new}.  The proof is standard and the details are left to the reader. We  note that the condition $t>1/2$ ensures the existence of the trace of $g_1$ on the boundary. \proofend

\subsection{Proof of Corollary~\ref{cor3}}
The proof is similar to the one of Corollary~\ref{cor1}. The details are left to the reader.  \proofend

\appendix

\section{Some averaging estimates}\label{apA}
\renewcommand{\theequation}{A\arabic{equation}}
\renewcommand{\thelemma}{A\arabic{lemma}}
  \setcounter{lemma}{0}  



\subsection{Proof of Lemma~\ref{lem-average1}} 

It is clear that 
\begin{equation*}
|\hat K_{s\sigma_{2} - i s \sigma_{1}} (k)|^{p} \le \frac{C_p}{|k|^{2p}} \mbox{ for } |k| > 2s. 
\end{equation*}
Hence to obtain \eqref{est1-lem-average1},  it suffices to prove that
\begin{equation*}
 \int_{\sigma_{1} \in \mS^{d-1}} \int_{\sigma_{2} \in \mS_{\sigma_{1}}} |\hat K_{s\sigma_{2} - i s \sigma_{1}} (k) |^{p} \, d \sigma_{2} \, d \sigma_{1} \le \frac{C_p}{|k|^{p} s^{p}} \quad \mbox{ for } |k| \le 2 s. 
\end{equation*}
Without loss of generality one may assume that  $k = t e_{1} \big( e_{1} =   (1, 0,  \cdots, 0 ) \big)$. Set 
\begin{equation*}
\xi = s (\sigma_{2} - i \sigma_{1}). 
\end{equation*}
We have
\begin{equation*}
 \frac{1}{\big|-|k|^{2} + i k \cdot \xi \big|^{p}}  =  \frac{1}{\big|-t^{2} + i s t e_{1} \cdot \sigma_{1} +  s t e_{1} \cdot \sigma_{2} \big|^{p}} \sim \frac{1}{|t^{2} - s t \sigma_{1} \cdot e_{1}|^{p} + (s t)^{p} |\sigma_{2} \cdot e_{1}|^{p}} .  
\end{equation*}
Let $\theta_{1}$ be the angle between $\sigma_{1}$ and $e_{1}$ and let $\theta_{2}$ be the angle between $\sigma_{2}$ and $v$ where  $v = e_1 - (e_1 \cdot \sigma_1) \sigma_1 = e_1 - \cos \theta_{1} \sigma_{1}$. Note that $v \in \mbox{span} \{ \sigma_{1}, e_{1} \}$, $v$ is orthogonal to $\sigma_{1}$, and $|v| = |\sin(\theta_1)|$. Using  the spherical area element, we have
\begin{multline*}
C_p \int_{\sigma_{1} \in \mS^{d-1}} \int_{\sigma_{2} \in \mS^{d-1}_{\sigma_{1}}} \frac{1}{|t^{2} - s t \sigma_{1} \cdot e_{1}|^{p} + (s t)^{p} |\sigma_{2} \cdot e_{1}|^{p}}\, d \sigma_{2} \, d \sigma_{1}  \\[6pt]
\le \int_{0}^{\pi/2} \int_{0}^{\pi} \frac{\theta_{1}^{d-2} \theta_{2}^{d-3}}{|t^{2} - s t \cos \theta_{1}|^{p} + (s t)^{p} |\sin \theta_{1} \cos \theta_{2}|^{p}}  \, d \theta_{2} \, d \theta_{1} \\[6pt] + 
\int_{\pi/2}^{\pi} \int_{0}^{\pi} \frac{ (\pi-\theta_{1})^{d-2} \theta_{2}^{d-3}}{|t^{2} - s t \cos \theta_{1}|^{p} + (s t)^{p} |\sin \theta_{1} \cos \theta_{2}|^{p}}  \, d \theta_{2} \, d \theta_{1}. 
\end{multline*}
Here we use $|\sigma_2 \cdot e_1| = |\sigma_2 \cdot v| = |\sin \theta_1 \cos \theta_2|$. 
It follows that 
\begin{multline*}
\int_{\sigma_{1} \in \mS^{d-1}} \int_{\sigma_{2} \in \mS^{d-1}_{\sigma_{1}}} \frac{1}{\big|-|k|^{2} + i k \cdot \xi \big|^{p}} \, d \sigma_{2} \, d \sigma_{1} \\[6pt]
 \le  {C_p\over (st)^p} \int_{0}^{\pi/2} \int_{0}^{\pi} \frac{\theta_{1}^{d-2} \theta_{2}^{d-3}}{| {t \over s} -  \cos \theta_{1}|^{p} +  |\sin \theta_{1} \cos \theta_{2}|^{p}}  \, d \theta_{2} \, d \theta_{1} . 
\end{multline*}
Fix $0< \delta < 2 - p$, and consider the case $t \leq s$.  Then
\begin{multline*}
 \int_{0}^{\pi/2} \int_{0}^{\pi} \frac{\theta_{1}^{d-2} \theta_{2}^{d-3}}{| {t \over s} -  \cos \theta_{1}|^{p} +  |\sin \theta_{1} \cos \theta_{2}|^{p}}  \, d \theta_{2} \, d \theta_{1} \\[6pt] 
\le C_p \int_{0}^{\pi/2}\int_0^{\pi}  \frac{\theta_{1}^{d-2} \theta_2^{d-3} }{|{t\over s} -  \cos \theta_{1}|^{1 - \delta} | \sin \theta_{1} \cos \theta_2 |^{p-1 + \delta}} \, d \theta_{2} \, d \theta_1 .
\end{multline*}
This implies 
\begin{multline*}
\int_{0}^{\pi/2} \int_{0}^{\pi} \frac{\theta_{1}^{d-2} \theta_{2}^{d-3}}{| {t \over s} -  \cos \theta_{1}|^{p} +  |\sin \theta_{1} \cos \theta_{2}|^{p}}  \, d \theta_{2} \, d \theta_{1}  \\[6pt] 
\le \int_{0}^{\pi/2}  \frac{\theta_{1}^{d-2}  }{|{t\over s} - \cos \theta_{1}|^{1-\delta} | \sin \theta_{1} |^{p-1+\delta}} \, d \theta_{1} \int_0^{\pi} {1 \over |\cos \theta_2|^{p-1+\delta}} d \theta_2
\end{multline*}
A computation yields 
\begin{equation}\label{m1}
\int_{0}^{\pi/2} \int_{0}^{\pi} \frac{\theta_{1}^{d-2} \theta_{2}^{d-3}}{|t^{2} - s t \cos \theta_{1}|^{p} + (s t)^{p} |\sin \theta_{1} \cos \theta_{2}|^{p}}  \, d \theta_{2} \, d \theta_{1}  \\[6pt] 
\le   C_p   \int_{0}^{\pi/2}  \frac{\theta_{1}^{d-1 - p - \delta}  }
{| {t \over s} -  \cos \theta_{1}|^{1-\delta}}  \, d \theta_{1} . 
\end{equation}
On the other hand,  let $\theta_0,\alpha_0$ be such that $\cos \theta_0 = {t \over s}$ and $|\cos \theta_0- \cos(\alpha + \theta_0) | \leq {1\over 2 }$ for all $| \alpha| \leq \alpha_0$. We have, since $d-1 -p - \delta \ge 2 - p - \delta > 0$,  
\begin{multline}\label{m2}
C_p \int_{0}^{\pi/2}  \frac{\theta_{1}^{d-1 - p - \delta}  }
{| {t \over s} -  \cos \theta_{1}|^{1-\delta}}  \, d \theta_{1} \\[6pt]
 \le  \int_{ | \theta - \theta_0| \leq \alpha_0 }  \frac{1  }
{| \cos \theta_0 -  \cos \theta|^{1-\delta}}  \, d \theta 
+   \int_{ [0,\pi/2] \backslash \{ | \theta - \theta_0| \leq \alpha_0\}  }  \frac{1  }
{| \cos \theta_0 -  \cos \theta|^{1-\delta}}  \, d \theta. 
\end{multline}
We have 
\begin{align}\label{m3}
\int_{ | \theta - \theta_0| \leq \alpha_0 }  \frac{1  }
{| \cos \theta_0 -  \cos \theta|^{1-\delta}}  \, d \theta &   \le   \int_{ | \theta - \theta_0| \leq \alpha_0 }  \frac{C_p  }
{|\sin [(\theta_0 + \theta)/2]|^{1- \delta}|\theta - \theta_0|^{1-\delta}}  \, d \theta  \nonumber \\[6pt]
& \le  \int_{ | \theta - \theta_0| \leq \alpha_0 }  \frac{C_p  }
{ |\sin \theta_0|^{1-\delta}|\theta - \theta_0|^{1-\delta}}  \, d \theta \le  \frac{C_p}{ (1 - ({t\over s}))^{1-\delta\over2} }
\end{align}
and 
\begin{equation}\label{m4}
  \int_{ [0,\pi/2] \backslash \{ | \theta - \theta_0| \leq \alpha_0\}  }  \frac{1  }
{| {t \over s} -  \cos \theta|^{1-\delta}}  \, d \theta \le C_p.
\end{equation}
A combination of \eqref{m1}, \eqref{m2}, \eqref{m3}, and \eqref{m4} yields 
\begin{equation}\label{aver.shell.est}
 \int_{0}^{\pi/2} \int_{0}^{\pi} \frac{\theta_{1}^{d-2} \theta_{2}^{d-3}}{|t^{2} - s t \cos \theta_{1}|^{p} + (s t)^{p} |\sin \theta_{1} \cos \theta_{2}|^{p}}  \, d \theta_{2} \, d \theta_{1}  \le \frac{C_p}{ (s t)^{p} (1 - ({t\over s}))^{1-\delta\over2} } + \frac{C_p}{ (s t)^{p}} . 
\end{equation}
For $s < t \le 2s$, we have 
\begin{equation*}
  \int_{0}^{\pi/2}  \frac{\theta_{1}^{d-1 - p - \delta}  }
{| {t \over s} -  \cos \theta_{1}|^{1-\delta}}  \, d \theta_{1}  \le C. 
\end{equation*}
Hence we also obtain \eqref{aver.shell.est} in this case.  Averaging \eqref{aver.shell.est} in $s$ yields  bound  \eqref{est1-lem-average1}.

\medskip We now  establish \eqref{est2-lem-average1}. Define $v_1 = e_{1} - (e_{1} \cdot \sigma_{1}) \sigma_{1} - (e_{1} \cdot \sigma_{2})  \sigma_{2} = v - (v \cdot \sigma_2) \sigma_2$ and let $\theta_{3}$ be the angle between $\sigma_3$ and $v_{1}$. We have, since $\sigma_3 \cdot e_1 = \sigma_3 \cdot v_1 = |v_1| \cos \theta_3$,  
\begin{multline*}
  \int_{\sigma_{1} \in \mS^{d-1}} \int_{\sigma_{2} \in \mS^{d-1}_{\sigma_{1}}} \int_{\sigma_{3} \in \mS^{d-1}_{\sigma_{1}, \sigma_{2}}} \Big|\hat K_{\frac{s^2\sigma_{2}}{ \sqrt{1 + s^{2}}} + \frac{s\sigma_{3}}{ \sqrt{1 + s^{2}} }  - i s \sigma_{1} } (k) \Big|^{p}  \, d \sigma_{3}\, d \sigma_{2} \, d \sigma_{1}\\[6pt]
\le C_p \int_{0}^{3\pi/4} \int_{0}^{\pi}  \int_{0}^{\pi} \frac{\theta_{1}^{d-2} \theta_{2}^{d-3} \theta_{3}^{d-4}}{|t^{2} - s t \cos \theta_{1}|^{p} + (st)^{p} |\sin \theta_{1} \cos \theta_{2} - |v_1| \cos \theta_{3} / s|^{p}} \, d \theta_{3} \, d \theta_{2} \, d \theta_{1}. 
\end{multline*}
Here $\int_{0}^{\pi} f(\theta_{3}) \theta_{3}^{d-4} \,  d \theta_{3} := f(\pi) + f (0)$ if $d=3$. We will only  consider the case $d \ge 4$, the case $d =3$ follows similarly. We have 
\begin{equation*}
  \int_{0}^{\pi} \frac{ \theta_{3}^{d-4} \, d \theta_{3}}{|t^{2} - s t \cos \theta_{1}|^{p} + (st)^{p} |\sin \theta_{1} \cos \theta_{2} - |v_1| \cos \theta_{3} / s|^{p}} \le \frac{C_p}{|t^{2} - s t \cos \theta_{1}|^{p-1} t |v_1|} .
\end{equation*}
Since 
\begin{equation*}
|v_1|^2 = |v|^2 - |v \cdot \sigma_2|^2 = \sin^2 \theta_1 \sin^2 \theta_2, 
\end{equation*}
it follows that 
\begin{equation*}
  \int_{0}^{\pi} \frac{ \theta_{3}^{d-4} \, d \theta_{3}}{|t^{2} - s t \cos \theta_{1}|^{p} + (st)^{p} |\sin \theta_{1} \cos \theta_{2} - |v_1| \cos \theta_{3} / s|^{p}} \le \frac{C_p}{|t^{2} - s t \cos \theta_{1}|^{p-1} t |\sin \theta_2| |\sin \theta_1|}. 
\end{equation*}
This implies 
\begin{multline*}
\int_{0}^{3\pi/4} \int_{0}^{\pi}  \int_{0}^{\pi} \frac{\theta_{1}^{d-2} \theta_{2}^{d-3} \theta_{3}^{d-4}}{|t^{2} - s t \cos \theta_{1}|^{p} + (st)^{p} |\sin \theta_{1} \cos \theta_{2} - |v_1| \cos \theta_{3} / s|^{p}} \, d \theta_{3} \, d \theta_{2} \, d \theta_{1}
\\[6pt]
 \le 
C_p  \int_{0}^{3\pi/4} \int_{0}^{\pi} \frac{\theta_{1}^{d-2} \theta_{2}^{d-3}}{|t^{2} - s t \cos \theta_{1}|^{p-1} t |\sin \theta_1 \sin \theta_2|}  \, d \theta_2 \, d \theta_1. 
\end{multline*}
We have,  since $d \ge 4$, 
\begin{equation*}
  \int_{0}^{3\pi/4} \int_{0}^{\pi} \frac{\theta_{1}^{d-2} \theta_{2}^{d-3}}{|t^{2} - s t \cos \theta_{1}|^{p-1} t |\sin \theta_1 \sin \theta_2|}  \, d \theta_2 \, d \theta_1 
\le \int_{0}^{3\pi/4}  \frac{1}{|t^{2} - s t \cos \theta_{1}|^{p-1} t } \, d \theta_1 
\le  \frac{C_p}{t^p s^p}.  
\end{equation*}
We obtain the conclusion. \proofend

\subsection{Proof of Lemma~\ref{lem-average-11}} We first claim that, for $k \in \mR^3$ with $|k| \ge 2$,  
 \begin{equation}\label{T0}
{1 \over R} \int_{R/2}^{2R} \mathop{ \int_{\sigma_{1} \in \mS^{2}} \int_{\sigma_{2} \in \mS^{2}_{\sigma_{1}}}}_{ \dist (k, \Gamma_{\xi}) \ge  |k| / |\xi|}  |\hat K_{s\sigma_{2} - i s \sigma_{1}} (k)|^{2} \, d \sigma_{2} \, d \sigma_{1}  \, ds \le C \min \Big\{ \frac{\ln R}{R^{2}|k|^{2}}, \frac{1}{|k|^{4}} \Big\}. 
 \end{equation}
 Here $\xi = \xi(s, \sigma_1, \sigma_2) = s \sigma_{2} - i s \sigma_{1}$ and 
 \begin{equation*}
 \Gamma_{\xi} := \{k \in \mR^{3} ; \; - |k|^{2}  + i \xi \cdot k =0 \} .
 \end{equation*}
Indeed, since
\begin{equation*}
|\hat K_{s\sigma_{2} - i s \sigma_{1}} (k)|^{2} \le \frac{C}{|k|^{4}} \mbox{ for } |k| > 2s,  
\end{equation*}
it suffices to prove that
\begin{equation}\label{T1}
{1\over R} \int_{R/2}^{2R} \mathop{\int_{\sigma_{1} \in \mS^{2}} \int_{\sigma_{2} \in \mS^2_{\sigma_{1}}}}_{ \dist (k, \Gamma) \ge  |k| / |\xi|}  |\hat K_{s\sigma_{2} - i s \sigma_{1}} (k) )|^{2} \, d \sigma_{2} \, d \sigma_{1} \, ds \le \frac{C \ln R}{|k|^{2} R^{2}} \quad \mbox{ for } |k| \le 2 s. 
\end{equation}
Without loss of generality, one may assume that  $k = t e_{1} =   (t, 0, 0 ) $.  As in the proof of Lemma~\ref{lem-average1}, we have 
\begin{multline}\label{T2}
\int_{\sigma_{1} \in \mS^{d-1}} \int_{\sigma_{2} \in \mS^{d-1}_{\sigma_{1}}} \frac{1}{\big|-|k|^{2} + i k \cdot \xi \big|^{2} + (t/s)^2} \, d \sigma_{2} \, d \sigma_{1} \\[6pt]
 \le  {C\over (st)^2} \int_{0}^{\pi/2} \int_{0}^{\pi} \frac{\theta_{1}}{| {t \over s} -  \cos \theta_{1}|^{2} +  |\sin \theta_{1} \cos \theta_{2}|^{2}  + s^{-4}}  \, d \theta_{2} \, d \theta_{1} . 
\end{multline}
A computation yields
\begin{equation}\label{T3}
 \int_{0}^{\pi/2} \int_{0}^{\pi} \frac{\theta_{1}}{| {t \over s} -  \cos \theta_{1}|^{2} +  |\sin \theta_{1} \cos \theta_{2}|^{2}  + s^{-4}}  \, d \theta_{2} \, d \theta_{1} \le C  \int_{0}^{\pi/2} \frac{1}{\big| \frac{t}{s} - \cos \theta_1\big| + s^{-2}} \, d \theta_1. 
\end{equation}
and 
\begin{equation}\label{T4}
\int_{0}^{\pi/2} \frac{1}{\big| \frac{t}{s} - \cos \theta_1\big| + s^{-2}} \, d \theta_1 \le C \ln s. 
\end{equation}
A combination of \eqref{T2}, \eqref{T3}, and \eqref{T4} yields \eqref{T1}; hence \eqref{T0} is established. 
\medskip

In the rest, we only give the proof of \eqref{avrg}. The proof of \eqref{avrg2} follows similarly. Applying \eqref{T0},  we have
\begin{multline*}
{1\over R} \int_{R/2}^{2R} 
\int_{\mS^{2}} \int_{\mS^{2}_{\sigma_{1}}} \int_{\mR^{3}} |\hat q(\eta)|^{2} \int_{4 {|\xi|} \ge \dist(k, \Gamma_{s \sigma_{2} - i s \sigma_{1}}) \ge |k|/ {|\xi|}} \frac{|k|^{2} | \hat K_{s \sigma_{2} - i s \sigma_{1}}(k)|^{2}}{|k - \eta|^{2}} \, dk \, d \eta \, d \sigma_2 \, d \sigma_1 \, d s \\[6pt]
\le  C \int_{\mR^{3}} |\hat q(\eta)|^{2} \int_{10 R \ge |k| } \frac{ \ln R}{R^{2 }|k - \eta|^{2}} \, dk \, d \eta. 
\end{multline*}
Since
\begin{equation*}
\int_{\mR^{3}} |\hat q(\eta)|^{2} \int_{10 R \ge |k| } \frac{ \ln R}{R^{2 }|k - \eta|^{2}} \, dk \, d \eta 
\le C \int_{\mR^{3}} |\hat q(\eta)|^{2} \min \Big\{\frac{\ln R}{R}, \frac{R\ln R}{|\eta|^{2}} \Big\} \, d \eta, 
\end{equation*}
it follows that 
\begin{multline}\label{point1}
{1\over R} \int_{R/2}^{2R} \int_{\mS^{2}} \int_{\mS^{2}_{\sigma_{1}}} \int_{\mR^{3}} |\hat q(\eta)|^{2} \int_{4 {|\xi|} \ge \dist(k, \Gamma_{s \sigma_{2} - i s \sigma_{1}}) \ge |k|/ {|\xi|}} \frac{|k|^{2} | \hat K_{s \sigma_{2} - i s \sigma_{1}}(k)|^{2}}{|k - \eta|^{2}} \, dk \, d \eta \, d \sigma_2 \, d \sigma_1 \, d s \\[6pt]
\le C \int_{\mR^{3}} |\hat q(\eta)|^{2}   \min \Big\{\frac{\ln R}{R}, \frac{R \ln R}{|\eta|^{2}} \Big\} \, d \eta. 
\end{multline}

Define 
\begin{equation*}
\widetilde q(k) : = \sup_{\eta \in B_{4}(k)} |\hat q(\eta)|. 
\end{equation*}
We have 
\begin{multline}\label{step1}
{1\over R} \int_{R/2}^{2R} \int_{\mS^{2}} \int_{\mS^{2}_{\sigma_{1}}}  \int_{\mR^{3}} |\hat q(\gamma)|^{2} \int_{\dist(\eta, \Gamma_{s \sigma_{2} - i s \sigma_{1}}) \le |\eta|/ {|\xi|}} \frac{1}{|\eta - \gamma|^{2}} \, d \eta \,  d \gamma \, d \sigma_2 \, d \sigma_1 \, ds \\[6pt]
\le C \int_{\mR^{3}} |\tilde q(\gamma)|^{2}  \min\Big\{\frac{R}{ |\gamma|^{2}}, \frac{1}{R} \Big\} \, d \gamma.
\end{multline}
Fix $q \in C^{\infty}_{\mc}(\mR^{3})$  such that $\varphi = 1$ in $B_{1}$ and $\supp \varphi \subset B_{2}$ and define 
\begin{equation*}
\widetilde \varphi (k) = \sup_{\eta \in B_{4}(k)} |\hat \varphi (\eta) |. 
\end{equation*}
Using the fact that 
\begin{equation}\label{dphi}
|\tilde q| \le \tilde \varphi * |\hat q|, 
\end{equation}
and applying H\"older's inequality, we have
\begin{equation*}
\int_{\mR^{3}} |\tilde q(\gamma)|^{2}\min\Big\{\frac{R}{ |\gamma|^{2}}, \frac{1}{R} \Big\} \, d \gamma \le C \int_{\mR^{3}}|\hat q(\beta)|^{2} \int_{\mR^{3}} \tilde \varphi (\gamma - \beta)\min\Big\{\frac{R}{ |\gamma|^{2}}, \frac{1}{R} \Big\} \, d \gamma \, d \beta. 
\end{equation*}
It follows from \eqref{dphi} that 
\begin{equation}\label{step2}
\int_{\mR^{3}} |\tilde q(\gamma)|^{2} \min\Big\{\frac{R}{ |\gamma|^{2}}, \frac{1}{R} \Big\} \, d \gamma \le C \int_{\mR^{3}} |\hat q(\beta)|^{2}  \min\Big\{\frac{R}{ |\beta|^{2}}, \frac{1}{R} \Big\} \,  d \beta. 
\end{equation}
A combination of \eqref{step1} and \eqref{step2} yields 
\begin{multline}\label{point2}
{1\over R} \int_{R/2}^{2R} \int_{\mS^{2}} \int_{\mS^{2}_{\sigma_{1}}}  \int_{\mR^{3}} |\tilde q(\gamma)|^{2} \int_{\dist(\eta, \Gamma_{s \sigma_{2} - i s \sigma_{1}}) \le |\eta|/ {|\xi|}} \frac{1}{|\eta - \gamma|^{2}} \, d \eta \,  d \gamma \, d \sigma_2 \, d \sigma_1 \, ds \\[6pt]
\le 
C \int_{\mR^{3}} |\hat q(\beta)|^{2}  \min\Big\{\frac{R}{ |\beta|^{2}}, \frac{1}{R} \Big\} \,  d \beta. 
\end{multline}
We derive \eqref{avrg} from \eqref{point1},  and \eqref{point2}.  \proofend

\section{Boundary determination}
\label{apB}
\renewcommand{\theequation}{B\arabic{equation}}
\renewcommand{\thelemma}{B\arabic{lemma}}
\renewcommand{\theproposition}{B\arabic{proposition}}

  \setcounter{lemma}{0}  
    \setcounter{proposition}{0}

In this appendix, we prove the following result

\begin{proposition}\label{Pro-boundary} Let $d \ge 2$, $\Omega$ be an open subset of $\mR^d$ of class $C^{1}$, and $\gamma_1, \gamma_2 \in W^{1,1}(\Omega)$. Assume $DtN_{\gamma_1} = DtN_{\gamma_2}$, then we have
\begin{equation*}
\gamma_1 = \gamma_2 \mbox{ on } \partial \Omega.
\end{equation*}
\end{proposition}

\noindent{\bf Proof.} We give the proof in the case $d \ge 3$. The proof in the $2d$ case follows similarly. We prove this result by contradiction. Assume that the conclusion is not true. Hence there exists some $z$ on $\partial \Omega$ such that
\begin{equation}\label{contrad-1}
\gamma_1(z) \neq \gamma_2(z)
\end{equation}
\begin{equation}\label{contrad-2}
\lim_{r \to 0} \mint_{B(z,r) \cap \Omega} |\gamma_1(x) -\gamma_1(z)| = 0,
\end{equation}
and
\begin{equation}\label{contrad-3}
\lim_{r \to 0} \mint_{B(z,r) \cap \Omega} |\gamma_2(x) -\gamma_2(z)| = 0.
\end{equation}
These last two statement following from the fact that for ${\mathcal H}^{d-2}$ a.e. $y \in \partial \Omega$, we have (see e.g. \cite[Theorem 2 on page 181]{EGMeasure})
\begin{equation*}
\lim_{r \to 0} \mint_{B(y,r) \cap \Omega} |\gamma_1(x) -\gamma_1(y)| = 0,
\end{equation*}
and
\begin{equation*}
\lim_{r \to 0} \mint_{B(y,r) \cap \Omega} |\gamma_2(x) -\gamma_2(y)| = 0.
\end{equation*}
Let $z_n$ be a sequence in $\mR^d \setminus \Omega$ such that
\begin{equation*}
\dist(z_n, \Omega) = |z_n - z| \quad \mbox{ and } \quad \lim_{n \to \infty} |z_n - z| =0.
\end{equation*}
Set
\begin{equation*}
v_n = \frac{1}{|x - z_n|^{d -2}} \mbox{ in } \mR^d,
\end{equation*}
and let $u_{j, n} \in H^1(\Omega)$ ($j=1,2$) be the unique solution to the system
\begin{equation*}\left\{
\begin{array}{cc}
\dive(\gamma_j \nabla u_{j, n}) = 0 & \mbox{ in } \Omega, \\[6pt]
u_{j, n} = v_n & \mbox{ on } \partial \Omega.
\end{array} \right.
\end{equation*}
Define
\begin{equation*}
w_{j, n} = u_{j, n} - v_n \mbox{ in } \Omega.
\end{equation*}
It is clear that
\begin{equation}\label{eq-vn}
\Delta v_n = 0 \quad \mbox{ in } \Omega.
\end{equation}
We also have
\begin{equation*}
- \dive(\gamma_j \nabla w_{j, n}) = - \dive(\gamma_j \nabla u_{j, n}) -  \dive(\gamma_j \nabla v_n) = - \dive([\gamma_j - \gamma_j(z)] \nabla v_n) \mbox{ in } \Omega,  
\end{equation*}
where in the last identity, we used \eqref{eq-vn}. This implies
\begin{equation*}
\int_{\Omega} \gamma_j |\nabla w_{j, n}|^2 = \int_{\Omega} [\gamma_j - \gamma_j(z)] \nabla v_n \nabla w_{j, n}.
\end{equation*}
It follows from \eqref{contrad-2} and \eqref{contrad-3} that
\begin{equation*}
\| \nabla w_{j, n}\|_{L^2} \le \| [\gamma_j - \gamma_j(z)] \nabla v_n\|_{L^2} = \frac{o(1)}{|z - z_n|^{(d-2)/2}}.
\end{equation*}
Here and in the following we let $o(1)$ denote a quantity going to 0 as $n \to \infty$; hence,
\begin{equation*}
\nabla w_{j, n} = \nabla v_n + \frac{g_n}{|z - z_n|^{(d-2)/2}},
\end{equation*}
for some  $\| g_n\|_{L^2} \to 0$ as $n \to \infty$. On the other hand,
\begin{equation*}
\int_{\Omega} (\gamma_{1} - \gamma_2) \nabla w_{1, n} \nabla w_{2, n} = 0
\end{equation*}
which implies
\begin{equation*}
[\gamma_1(z) - \gamma_2(z)] \frac{1}{|z-z_n|^{d-2}} = o(1)\frac{1}{|z - z_n|^{d-2}}.
\end{equation*}
Hence
\begin{equation*}
\gamma_1(z) = \gamma_2(z).
\end{equation*}
This contradicts \eqref{contrad-1}, and the conclusion follows. \proofend
\vskip1cm

\noindent {\bf{Acknowledgements}}
{
Hoai-Minh Nguyen was supported in part by NSF grant DMS-1201370 and by the Alfred P. Sloan Foundation.
Daniel Spirn was supported in part by NSF grant DMS-0955687.  We would like to  thank Jean-Pierre Puel for pointing out an error in the proof of Theorem 2 in an earlier version. 
}

\end{document}